\documentclass[12pt,a4paper]{article}
\usepackage{mathrsfs}
\usepackage{amssymb}
\usepackage{amsmath}
\usepackage{amsfonts}
\usepackage{amstext}
\usepackage{amsthm}
\usepackage{graphicx}
\usepackage[top=1in, bottom=1in,left=0.75in,right=0.75in]{geometry}

\def \qed {\hfill \vrule height6pt width 6pt depth 0pt}

\begin{document}

\title{Eigenvalue, bifurcation, existence and nonexistence of solutions for Monge-Amp\`{e}re equations
\thanks{Research supported by the NSFC (No. 11061030).}}
\author{{\small  Guowei Dai\thanks{Corresponding author. Tel: +86 931
7971124.\newline
\text{\quad\,\, E-mail address}: daiguowei@nwnu.edu.cn.
}
} \\
{\small Department of Mathematics, Northwest Normal
University, Lanzhou, 730070, P.R. China}\\
}
\date{}
\maketitle

\begin{abstract}
In this paper we study the following eigenvalue boundary value problem for Monge-Amp\`{e}re equations:
\begin{equation}
\left\{
\begin{array}{l}
\det\left(D^2u\right)=\lambda^N f(-u)\,\, \text{in}\,\, \Omega,\\
u=0~~~~~~~~~~~~~~~~~~~~~~\,\,\text{on}\,\, \partial \Omega.
\end{array}
\right.\nonumber
\end{equation}
We establish the unilateral global bifurcation results for the problem with $f(u)=u^N+g(u)$ and $\Omega$ being the unit ball of $\mathbb{R}^N$.
More precisely, under some natural hypotheses on the perturbation function $g:\mathbb{R}\rightarrow\mathbb{R}$, we show that
$\left(\lambda_1,0\right)$ is a bifurcation point of the problem and there are two distinct unbounded continua of one-sign solutions,
where $\lambda_1$ is the first eigenvalue of the problem with $f(u)=u^N$. As the applications of the above results, we consider
with determining interval of $\lambda$,
in which there exist solutions for this problem in unit ball.
Moreover, we also get some results on the existence and nonexistence of convex solutions for this problem in general domain by domain comparison method.
\end{abstract}\textbf{\ }

\numberwithin{equation}{section}

\numberwithin{equation}{section}

\section{Introduction}

\quad\, The Monge-Amp\`{e}re equations are a type of important fully nonlinear elliptic equations [\ref{GT}, \ref{T}]. The study
of Monge-Amp\`{e}re equations has been received considerable attention in recent years. Historically, the
study of Monge-Amp\`{e}re equations is motivated by Minkowski problem and Weyl problem.
Existence and regularity results may be found in [\ref{CNS}, \ref{CY}, \ref{CY1}, \ref{GT}, \ref{K}, \ref{L},
\ref{P}, \ref{P1}, \ref{P2}, \ref{ZW}] and the reference therein.

We consider the following real Monge-Amp\`{e}re equations
\begin{equation}\label{MB}
\left\{
\begin{array}{l}
\det\left(D^2u\right)=\lambda^N f(-u)\,\, \text{in}\,\, B,\\
u=0~~~~~~~~~~~~~~~~~~~~~\,\,\text{on}\,\, \partial B,
\end{array}
\right.
\end{equation}
where $D^2u=\left(\frac{\partial u}{\partial x_i\partial x_j}\right)$ is the Hessian
matrix of $u$, $B$ is the unit ball of $\mathbb{R}^N$, $\lambda$ is a
nonnegative parameter and $f:\mathbb{R}\rightarrow \mathbb{R}$
is a continuous function. The study of problem (\ref{MB}) in
general domains of $\mathbb{R}^N$ may be found in [\ref{CNS}, \ref{GT}]. Kutev
[\ref{K1}] investigated the existence of strictly convex radial solutions
of problem (\ref{MB}) when $f(s)=s^p$. Delano [\ref{Del}] treated the
existence of convex radial solutions of problem (\ref{MB}) for a class of
more general functions, namely $\lambda\exp f(\vert x\vert,u,\vert \nabla u\vert)$.

In [\ref{HW}, \ref{K1}], the authors have showed that problem (\ref{MB})
can reduce to the following boundary value problem:
\begin{equation}\label{MO}
\left\{
\begin{array}{l}
\left(\left(u'\right)^N\right)'=\lambda^NNr^{N-1} f(-u)\,\, \text{in}\,\, 0<r<1,\\
u'(0)=u(1)=0.
\end{array}
\right.
\end{equation}
By a solution of problem (\ref{MO}) we understand it is a
function which belongs to $C^2[0,1]$ and satisfies (\ref{MO}). It has been known that any positive solution
of problem (\ref{MO}) is strictly concave and any negative solution
is strictly convex in $(0,1)$ so long as $f$ does not vanish on any entire interval (see [\ref{HW}]).
Under the assumption of $f\geq 0$, Wang [\ref{W}], Hu and Wang [\ref{HW}] also established several
criteria for the existence, multiplicity and nonexistence of strictly convex solutions for problem (\ref{MO}) using
fixed index theorem. However, there is no any information on the bifurcation points and the
optimal intervals for the parameter $\lambda$ so as to ensure existence of single or multiple solutions.
Fortunately, Lions [\ref{L1}] have proved the existence of the first eigenvalue $\lambda_1$ of problem (\ref{MB})
with $f(u)=u^N$ via constructive proof.

The first bifurcation phenomena in nonlinear problems is the bucking of the Euler rod, which proposed by Euler
in 1744. While the concept of bifurcation was firstly proposed by H. Poincar\'{e} in 1885. There are also
various concrete problems in the natural sciences involving of bifurcation phenomena, for example, Taylor
vortices [\ref{Bg}] and catastrophic shifts in ecosystems [\ref{S1}]. In this celebrated work [\ref{R}],
Rabinowitz established a unilateral global bifurcation theorem.
However, as pointed out by Dancer [\ref{D1}, \ref{D2}] and L\'{o}pez-G\'{o}mez [\ref{LG}], the proofs of these
theorems contain gaps. Fortunately, Dancer [\ref{D1}] gave a corrected version unilateral global bifurcation
theorem for linear operator which has been extended to the one-dimensional $p$-Laplacian problem by Dai and Ma [\ref{DM}].

Motivated by above, we shall establish a unilateral global bifurcation theorem
for problem (\ref{MO}) with $f(u)=u^N+g(u)$, i.e.,
\begin{equation}\label{Mg}
\left\{
\begin{array}{l}
\left(\left(u'\right)^N\right)'=\lambda^NNr^{N-1} \left(\left(-u\right)^N+g(-u)\right)\,\, \text{in}\,\, 0<r<1,\\
u'(0)=u(1)=0,
\end{array}
\right.
\end{equation}
where $g:\mathbb{R}\rightarrow \mathbb{R}$ satisfies $\lim_{ s\rightarrow0}g(s)/s^N=0$.
Concretely, we shall show that $\left(\lambda_1,0\right)$ is a bifurcation point of problem (\ref{Mg}) and there are two distinct
unbounded continua of one-sign solutions.

In global bifurcation theory of differential equations,
it is well known that a change of the index of the trivial solution implies
the existence of a branch of nontrivial solutions, bifurcating from the set of trivial solutions
and which is either unbounded or returns to the set of trivial solution.
Hence, the index formula of an isolated zero is very important in
the study of the bifurcation phenomena for semi-linear differential equations.
However, problem (\ref{Mg}) is a type of nonlinear equation. Hence, the common index formula
involving of linear map cannot be used here. In order to overcome this difficulty, we shall study an auxiliary
eigenvalue problem, which has an independent interesting, and establish an
index formula for it. Then by use of the index formula about of the auxiliary problem, we prove an index formula
involving of problem (\ref{Mg}) which guarantees $\left(\lambda_1,0\right)$ is a bifurcation point of nontrivial solutions
to problem (\ref{Mg}). Furthermore, by the similar arguments to the proofs of [\ref{DM}], we can get unilateral
global bifurcation results for problem (\ref{Mg}).

Based on the above unilateral global bifurcation results, we investigate the existence of strictly convex
or concave solutions of problem (\ref{MO}). We shall give the optimal intervals for the parameter $\lambda$
so as to ensure existence of single or multiple strictly convex or concave solutions. In order to study the exact
multiplicity of one-sign solutions for problem (\ref{MO}), we introduce the concept of stable solution.
Then by Implicit Function Theorem and stability properties, under some more strict assumptions of $f$,
we can show that the nontrivial solutions branch of problem (\ref{MO}) can be a smooth curve. Our results extend
the corresponding results to [\ref{HW}, \ref{L1}, \ref{W}].

On the basis of results on unit ball, we also study problem (\ref{MB}) on a general domain $\Omega$, i.e.,
\begin{equation}\label{MAh1}
\left\{
\begin{array}{l}
\det\left(D^2u\right)=\lambda^N f(-u)\,\, \text{in}\,\, \Omega,\\
u=0~~~~~~~~~~~~~~~~~~~~~\,\,\text{on}\,\, \partial \Omega,
\end{array}
\right.
\end{equation}
where $\Omega$ is a bounded convex domain of $\mathbb{R}^N$ with smooth boundary and $0\in\text{Int}\Omega$.
It is well-known [\ref{GT}] that problem (\ref{MAh1}) is elliptic only when the Hessian matrix $D^2u$ is positive (or negative) definite and
it is therefore natural to confine our attention to convex (or concave) solutions and positive (or negative) functions $f$.
Obviously, any convex solution of problem (\ref{MAh1}) is negative and strictly convex. In [\ref{ZW}], the authors has proved a lemma concerning the comparison between domains for problem (\ref{MAh1}) with
$f(s)=e^s$ by sub-supersolution method. We shall show that this lemma is also valid for problem (\ref{MAh1}).
Using this domain comparison lemma and the results on unit ball, we can prove some existence and nonexistence of
solutions for problem (\ref{MAh1}).

The rest of this paper is arranged as follows. In Section 2, we study an auxiliary problem and prove a
key index formula. In Section 3, we establish a unilateral global
bifurcation theorem for problem (\ref{Mg}). In Section 4, we give the intervals for the parameter $\lambda$
which ensure existence of single or multiple strictly convex or concave solutions for problem (\ref{MO}) under
some suitable assumptions of nonlinearity $f$. In Section 5, under some more strict assumptions of $f$,
we prove the exact multiplicity of one-sign solutions for problem (\ref{MO}). In Section 6, we prove some existence and nonexistence of convex
solutions for problem (\ref{MAh1}). \\

\section{A key preliminarily result}

\bigskip

\quad\, In this section, we shall study an auxiliary eigenvalue problem and prove a
key index formula that will be used in the next section.

Let $p\in[2,+\infty)$. Consider the following auxiliary problem
\begin{equation}\label{AP}
\left\{
\begin{array}{l}
-\left(\left\vert v'(r)\right\vert^{p-2}v'(r)\right)'=
\mu^{p-1} (p-1)r^{p-2}\vert v(r)\vert^{p-2}v(r)\,\, \text{in}\,\, 0<r<1,\\
v'(0)=v(1)=0.
\end{array}
\right.
\end{equation}
Let $X$ be the Banach space $C[0,1]$ with the norm
\begin{equation}
\Vert v\Vert=\sup_{r\in[0,1]}\vert v(r)\vert.\nonumber
\end{equation}
Define the map $T_\mu^p:X\rightarrow X$ by
\begin{equation}
T_\mu^pv=\int_1^r\varphi_{p'}\left(\int_s^0 \mu^{p-1}
(p-1)\tau^{p-2}\varphi_p(v)\,d\tau\right)\,ds,\,\, 0\leq r\leq1,\nonumber
\end{equation}
where $\varphi_p(s)=\vert s\vert^{p-2}s$, $p'=p/(p-1)$.
It is not difficult to verify that $T_\mu^p$ is continuous and compact.
Clearly, problem (\ref{AP}) can be equivalently written as
\begin{equation}
v=T_\mu^p v.\nonumber
\end{equation}

\indent Firstly, we show that the existence and uniqueness theorem is valid for problem (\ref{AP}).
\\ \\
\textbf{Lemma 2.1.} \emph{If $(\mu, v)$ is a solution of (\ref{AP})
and $v$ has a double zero, then $v \equiv 0$.}
\\ \\
\textbf{Proof.} Let $v$ be a solution of problem (\ref{AP}) and $r_*\in[0, 1]$ be a double zero.
We note that $v$ satisfies
\begin{equation}
v(r)=\int_{r_*}^r\varphi_{p'}\left(\int_s^{r_*} (p-1)\mu^{p-1}\tau^{p-2}\varphi_p(v)\,d\tau\right)\,ds.\nonumber
\end{equation}
Firstly, we consider $r\in[0, r_*]$. Then we have
\begin{eqnarray}
\vert v(r)\vert&\leq&\varphi_{p'}\left(\int_{r}^{r_*}(p-1)
\mu^{p-1}\tau^{p-2}\varphi_p(\vert v\vert)\,d\tau\right).\nonumber
\end{eqnarray}
Furthermore, it follows from above that
\begin{eqnarray}
\varphi_p(\vert v\vert)&\leq&\mu^{p-1}\int_{r}^{r_*}(p-1)\tau^{p-2}\varphi_p(\vert v\vert)\,d\tau. \nonumber
\end{eqnarray}
By the modification of Gronwall-Bellman inequality [\ref{ILL}, Lemma 2.2], we get $v \equiv 0$ on $[0, r^*]$.
Similarly, we can get $v \equiv 0$ on $[r^*, 1]$
and the proof is completed.\qed\\

Set $W_c^{1,p}(0,1):=\{v\in W^{1,p}(0,1)\big|v'(0)=v(1)=0\}$ with the norm
\begin{equation}
\Vert v\Vert_w=\left(\int_0^1\vert v'\vert^p\,dr\right)^{\frac{1}{p}}
+\left(\int_0^1(p-1)r^{p-2}\vert v\vert^p\,dr\right)^{\frac{1}{p}}.\nonumber
\end{equation}
Then it is easy to verify that $W_c^{1,p}(0,1)$ is a real Banach space.
\\ \\
\textbf{Definition 2.1.} We call that $v\in W_c^{1,p}(0,1)$ is the weak solution of
problem (\ref{AP}), if
\begin{equation}
\int_{0}^1\left\vert v'\right\vert^{p-2}v'\phi'\,{d}r=
(p-1)\mu^{p-1}\int_{0}^1r^{p-2}\vert v\vert^{p-2}v\phi\,{d}r\nonumber
\end{equation}
for any $\phi\in W_c^{1,p}(0,1)$.\\
\\
\indent For the regularity of weak solution, we have the following result.\\
\\
\textbf{Lemma 2.2.} \emph{Let $v$ be a weak solution of problem (\ref{AP}), then $v$ satisfies problem (\ref{AP}).}\\
\\
\indent In order to prove Lemma 2.2, we need the following technical result.
\\ \\
\noindent\textbf{Proposition 2.1.} \emph{Let $f:\mathbb{R}\rightarrow \mathbb{R}$ be a
function. For a given $x_0\in \mathbb{R}$, if $f$
is continuous in some neighborhood $U$ of $x_0$, differential in $U\setminus \{x_0\}$
and $\underset{x\rightarrow x_0}\lim f'(x)$ exists,
then $f$ is differential at $x_0$ and $f'(x_0)=\underset{x\rightarrow x_0}\lim f'(x)$.}\\ \\
\textbf{Proof.} The conclusion is a direct corollary of Lagrange mean Theorem, we omit its proof here.\qed
\\ \\
\textbf{Proof of Lemma 2.2.} According to Definition 2.1, we have
\begin{equation}
-\left(\left\vert v'(r)\right\vert^{p-2}v'(r)\right)'=
\mu^{p-1}(p-1)r^{p-2}\vert v(r)\vert^{p-2}v(r)\,\,\text{in}\,\,(0,1)\nonumber
\end{equation}
in the sense of distribution, i.e.,
\begin{equation}
-\left(\left\vert v'(r)\right\vert^{p-2}v'(r)\right)'
=\mu^{p-1}(p-1)r^{p-2}\vert v(r)\vert^{p-2}v(r)\,\,\text{in}\,\,(0,1)\setminus I\nonumber
\end{equation}
for some $I\subset (0,1)$ which satisfies $\text{meas}\{I\}=0$.
Furthermore, by virtue of the compact embedding of $W_c^{1,p}(0,1)\hookrightarrow
C^{\alpha}[0,1]$ with some $\alpha\in(0,1)$ (see [\ref{E}]), we obtain
\begin{equation}
-\left(\left\vert v'(r)\right\vert^{p-2}v'(r)\right)'\in C([0,1]\setminus I).\nonumber
\end{equation}
Let $u:=-\varphi_p\left(v'\right)$. The above relation follows that
$\lim_{r\rightarrow r_0} u'(r)$ exists for any $r_0\in I$. Thus,
Proposition 2.1 follows that $u\in C^1(0,1)$,
which implies that $v$ satisfies problem (\ref{AP}).\qed\\

\indent Define the functional $J$ on $W_c^{1,p}(0,1)$ by
\begin{equation}
J(v)=\int_0^1\frac{1}{p}\left\vert v'(r)\right\vert^p\,dr
-\mu^{p-1}\frac{p-1}{p} \int_0^1r^{p-2}\vert v\vert^p\,dr.\nonumber
\end{equation}
It is not difficult to verify that the critical points of $J$ are the weak solutions of problem (\ref{AP}).

Taking $f_1(v):=\int_0^1\frac{1}{p}\left\vert v'(r)\right\vert^p\,dr$ and
$f_2(v):=\frac{p-1}{p}\int_0^1r^{p-2}\vert v\vert^p\,dr$,
consider the following eigenvalue problem
\begin{equation}\label{APE}
A(v)=\eta Bv,
\end{equation}
where $A=\partial f_1$ and $B=\partial f_2$ denote the sub-differential of $f_1$ and $f_2$, respectively
(refer to [\ref{C}] for the details about of sub-differential).
By some simple computations, we can show that
\begin{equation}
\frac{f_1(v)}{f_2(v)}\geq \frac{1}{(p-1)}\nonumber
\end{equation}
for any $v\in W_c^{1,p}(0,1)$ and $v\not \equiv0$.
It is easy to see that the results of [\ref{IO}] remain true if (A0)$'$ is substituted by the following property:\\

(A0)$''$ \emph{Every positive solution of $(\ref{APE})$ satisfies $v\in C^1[0,1]$ and $v'(1)<0$.}\\

In view of Lemma 2.2, we can easily verify that all conditions of Theorem 1
and Theorem 2 of [\ref{IO}] and the assumption (A0)$''$ are satisfied.
Hence, by Theorem 1 and Theorem 2 of [\ref{IO}], we have the following result.
\\ \\
\textbf{Lemma 2.3.} \emph{Put $\eta_1(p)=\inf_{v\in W_c^{1,p}(0,1),v\not\equiv 0}\frac{f_1(v)}{f_2(v)}$.
Then we have} \\

\emph{(i) (\ref{APE}) has no nontrivial solution for $\eta\in(0,\eta_1(p))$,}

\emph{(ii) $\eta_1(p)$ is simple,}

\emph{(iii) (\ref{APE}) has a positive solution if and only if $\eta=\eta_1(p)$.}
\\ \\
\indent Let $\eta=\mu^{p-1}$, Lemma 2.3 shows the following result.\\
\\
\textbf{Lemma 2.4.} \emph{Put $\mu_1(p)=\left(\eta_1(p)\right)^{1/(p-1)}$.
Then we have} \\

\emph{(i) (\ref{AP}) has no nontrivial solution for $\eta\in(0,\mu_1(p))$,}

\emph{(ii) $\mu_1(p)$ is simple,}

\emph{(iii) (\ref{AP}) has a positive solution if and only if $\mu=\mu_1(p)$.}
\\ \\
\indent Moreover, we have the following result.
\\ \\
\textbf{Lemma 2.5.} \emph{If $(\mu,u)$ satisfies (\ref{AP}) and $\mu\neq \mu_1(p)$,
then $u$ must change sign.}\\
\\
\textbf{Proof.} Suppose that $u$ is not changing-sign. Without loss of
generality, we can assume that $u\geq 0$ in
$(0,1)$. Lemma 2.1 and 2.2 imply that $u>0$ in $(0,1)$. Lemma 2.4 implies $\mu=\mu_1(p)$ and $u=c v_1$
for some positive constant $c$, where $v_1$ is the positive eigenfunction
corresponding to $\mu_1(p)$ with $\left\Vert v_1\right\Vert=1$. This is a contradiction.\qed\\

\indent In addition, we also have that $\mu_1(p)$ is also isolated.
\\ \\
\textbf{Lemma 2.6.} \emph{$\mu_1(p)$ is the unique eigenvalue in $(0,\delta)$
for some $\delta>\lambda_1$.}
\\ \\
\textbf{Proof.} Lemma 2.4 has shown that $\mu_1(p)$ is left-isolated.
Assume by contradiction that there exists a sequence of eigenvalues
$\lambda_n\in(\mu_1(p), \delta)$ which converge to $\mu_1(p)$. Let $v_n$ be the
corresponding eigenfunctions.
Define
\begin{equation}
\psi_n:=\frac{v_n}{\left((p-1)\int_0^1 r^{p-2}\left\vert v_n\right\vert^p\,dr\right)^{\frac{1}{p}}}.\nonumber
\end{equation}
Clearly, $\psi_n$ are bounded in $W_c^{1,p}(0,1)$ so there exists a subsequence, denoted
again by $\psi_n$, and $\psi\in W_c^{1,p}(0,1)$ such that
$\psi_n\rightharpoonup \psi$ in $W_c^{1,p}(0,1)$ and $\psi_n\rightarrow \psi$ in $C^\alpha[0,1]$.
Since functional $f_1$ is sequentially weakly lower semi-continuous, we have
\begin{equation}
\int_0^1 \vert \psi'\vert^p\,dr \leq\liminf_{n\rightarrow+\infty}\int_0^1
\left\vert \psi_n'\right\vert^p\,dr=\mu_1^{p-1}(p).\nonumber
\end{equation}
On the other hand, $(p-1)\int_0^1 r^{p-2}\left\vert \psi_n\right\vert^p\,dr=1$ and
$\psi_n\rightarrow \psi$ in $C^\alpha[0,1]$ imply that
$(p-1)\int_0^1 r^{p-2}\vert \psi\vert^p\,dr=1$. Hence, $\int_0^1 \left\vert
\psi'\right\vert^p\,dr=\eta_1(p)$ via Lemma 2.3.
Then Lemma 2.1 and 2.3 show that $\psi>0$ in $(0,1)$. Thus $\psi_n\geq 0$
for $n$ large enough which contradicts the conclusion of Lemma 2.5.\qed
\\ \\
\indent Next, we show that the principle eigenvalue function $\mu_1:
[2,+\infty)\rightarrow \mathbb{R}$ is continuous.
\\ \\
\textbf{Lemma 2.7.} \emph{The eigenvalue function $\mu_1:[2,+\infty)\rightarrow \mathbb{R}$
is continuous.}
\\ \\
\textbf{Proof.} It is sufficient to show that $\eta_1(p):[2,+\infty)\rightarrow \mathbb{R}$
is continuous because of $\mu_1(p)=\left(\eta_1(p)\right)^{1/(p-1)}$.

From the variational characterization of $\eta_1(p)$ it follows that
\begin{equation}\label{cp1}
\eta_1(p)=\sup\left\{\lambda>0\Big|\lambda(p-1)\int_0^1r^{p-2}\vert v\vert^p
\,dr\leq\int_0^1\left\vert v'\right\vert^p\,dr
\,\,\text{for all\,\,}v\in C_c^\infty[0,1]\right\},
\end{equation}
where $C_c^\infty[0,1]=\left\{v\in C^\infty[0,1]\big|v'(0)=v(1)=0\right\}$, as $C_c^\infty[0,1]$
is dense in $W_c^{1,p}(0,1)$ (see [\ref{A}]).

Let $\left\{p_j\right\}_{j=1}^\infty$ be a sequence in $[2, +\infty)$ convergent to $p\geq2$. We shall show
that
\begin{equation}\label{cp2}
\lim_{j\rightarrow+\infty}\eta_1\left(p_j\right)=\eta_1(p).
\end{equation}
To do this, let $v\in C_c^\infty[0,1]$. Then, due to (\ref{cp1}), we get
\begin{equation}
\eta_1\left(p_j\right) \left(p_j-1\right)\int_0^1r^{p_j-2}\vert v\vert^{p_j}\,dr\leq\int_0^1\left\vert v'\right\vert^{p_j}\,dr.\nonumber
\end{equation}
On applying the Dominated Convergence Theorem we find
\begin{equation}\label{cp3}
\limsup_{j\rightarrow+\infty}\eta_1\left(p_j\right) (p-1)\int_0^1r^{p-2}\vert
v\vert^{p}\,dr\leq\int_0^1\left\vert v'\right\vert^{p}\,dr.
\end{equation}
Relation (\ref{cp3}), the fact that $v$ is arbitrary and (\ref{cp1}) yield
\begin{equation}
\limsup_{j\rightarrow+\infty}\eta_1\left(p_j\right)\leq\eta_1(p).\nonumber
\end{equation}

Thus, to prove (\ref{cp2}) it suffices to show that
\begin{equation}\label{cp4}
\liminf_{j\rightarrow+\infty}\eta_1\left(p_j\right)\geq\eta_1(p).
\end{equation}
Let $\left\{p_k\right\}_{k=1}^\infty$ be a subsequence of $\left\{p_j\right\}_{j=1}^\infty$ such that
$\underset{k\rightarrow+\infty}\lim\eta_1\left(p_k\right)=\underset{j\rightarrow+\infty}\liminf\eta_1\left(p_j\right)$.

Let us fix $\varepsilon_0>0$ so that $p-\varepsilon_0>1$ and for each
$0<\varepsilon<\varepsilon_0$ and $k\in \mathbb{N}$ large enough,
$p-\varepsilon<p_k<p+\varepsilon$. For $k\in \mathbb{N}$, let us choose $v_k\in W_c^{1,p_k}(0,1)$ such that $v_k>0$ in $(0,1)$,
\begin{equation}\label{cp5}
\int_0^1\left\vert v_k'\right\vert^{p_k}\,dr=1
\end{equation}
and
\begin{equation}\label{cp6}
\int_0^1\left\vert v_k'\right\vert^{p_k}\,dr=\eta_1\left(p_k\right)\left(p_k-1\right)\int_0^1 r^{p_k-2}\left\vert v_k\right\vert^{p_k}\,dr.
\end{equation}
For $0<\varepsilon<\varepsilon_0$ and $k\in \mathbb{N}$ large enough, (\ref{cp4}), (\ref{cp5})
and (\ref{cp6}) imply that
\begin{equation}\label{cp7}
\left\Vert v_k\right\Vert_{W_c^{1,p_k}(0,1)}\leq 1+\max\left\{\left(\frac{1}{\underset{k\rightarrow+\infty}
\lim\eta_1\left(p_k\right)}\right)^{\frac{1}{p+\varepsilon}},\left(\frac{1}
{\underset{k\rightarrow+\infty}\lim\eta_1\left(p_k\right)}\right)^{\frac{1}{p-\varepsilon}}\right\}.
\end{equation}
This shows that $\left\{v_k\right\}_{k=1}^\infty$ is a bounded sequence in $W_c^{1,p_k}(0,1)$, hence, in
$W_c^{1,p-\varepsilon}(0,1)$. Passing to a
subsequence if necessary, we can assume that $v_k \rightharpoonup v$ in
$W_c^{1,p-\varepsilon}(0,1)$ and hence
that $v_k \rightarrow v$ in $C^{\alpha}[0,1]$ with $\alpha=1-1/(p-\varepsilon)$
because the embedding of $W^{1,p-\varepsilon}(0,1)\hookrightarrow
C^{\alpha}[0,1]$ is compact. Thus,
\begin{equation}\label{cp8}
\left\vert v_k\right\vert^{p_k}\rightarrow\vert v\vert^{p}.
\end{equation}

We note that (\ref{cp6}) implies that
\begin{equation}\label{cp9}
\eta_1\left(p_k\right)\left(p_k-1\right)\int_0^1 r^{p_k-2}\left\vert v_k\right\vert^{p_k}\,dr=1
\end{equation}
for all $k\in \mathbb{N}$. Thus letting $k\rightarrow+\infty$ in (\ref{cp9}) and using (\ref{cp8}), we find
\begin{equation}\label{cp10}
\liminf_{j\rightarrow+\infty}\eta_1\left(p_j\right)(p-1)\int_0^1 r^{p-2}\vert v\vert^{p}\,dr=1.
\end{equation}
On the other hand, since $v_k\rightharpoonup v$ in $W_c^{1,p-\varepsilon}(0,1)$, from (\ref{cp5})
and H\"{o}lder's inequality we obtain that
\begin{equation}
\left\Vert v'\right\Vert_{p-\varepsilon}^{p-\varepsilon}\leq\liminf_{k\rightarrow+\infty}\left\Vert v_k'
\right\Vert_{p-\varepsilon}^{p-\varepsilon}\leq1,\nonumber
\end{equation}
where $\Vert \cdot\Vert_p$ denotes the normal of $L^p(0,1)$.
Now, letting $\varepsilon\rightarrow 0^+$, we find
\begin{equation}\label{cp11}
\left\Vert v'\right\Vert_p\leq1.
\end{equation}
Clearly, (\ref{cp10}), (\ref{cp11}) and $v\in W_c^{1,p-\varepsilon}(0,1)$ follow
that $v\in W_c^{1,p}(0,1)$.

Consequently, combining (\ref{cp10}) and (\ref{cp11}) we obtain
\begin{equation}
\liminf_{j\rightarrow+\infty}\eta_1\left(p_j\right)(p-1)\int_0^1 r^{p-2}\vert v\vert^{p}\,dr\geq
\int_0^1 \left\vert v'\right\vert^{p}\,dr.\nonumber
\end{equation}
This together with the variational characterization of $\eta_1(p)$ implies (\ref{cp4}) and hence
(\ref{cp2}). This concludes the proof of the lemma.\qed\\

We have known that $I-T^p_\mu$ is a completely continuous vector field in
$X$. Thus, the Leray-Schauder degree
$\deg\left(I-T^p_\mu, B_r(0),0\right)$ is well defined for
arbitrary $r$-ball $B_r(0)$ and $\mu\in (0,\delta)\setminus\left\{\mu_1(p)\right\}$, where $\delta$ comes from Lemma 2.6.
Now, we can compute it by the deformation along $p$.\\ \\
\noindent\textbf{Theorem 2.1.} \emph{Let $\mu$ be a constant with
$\mu\in(0,\delta)\setminus\left\{\mu_1(p)\right\}$. Then for
arbitrary $r>0$,}
\begin{equation}
\deg \left(I-T^p_\mu, B_r(0),0 \right)=\left\{
\begin{array}{l}
1, \,\,\,\,\,\,\text{if}\,\,\mu\in
\left(0,\mu_1(p)\right),\\
-1,\,\,\text{if}\,\,\mu\in\left(\mu_1(p),\delta\right).
\end{array}
\right.\nonumber
\end{equation}
\noindent\textbf{Proof.} We only treat the case of $\mu>\mu_1(p)$
because the proof for the case of $\mu<\mu_1(p)$ can be given similarly. Assume that $\mu_{1}(p) < \mu<\delta$.
Since the principle eigenvalue
depends continuously on $p$, there exist a continuous function
$\chi:[2,p]\rightarrow\mathbb{R}$ and $q\in [2,p]$ such that
$\mu_{1}(q) < \chi(q) <\delta$ and $\mu=\chi(p)$.
Define
\begin{equation}
\Upsilon(q,v)=v-\int_1^r\varphi_{p'}\left(\int_s^0 \left(\chi(q)\right)^{q-1}
(q-1)\tau^{q-2}\varphi_q(v)\,d\tau\right)\,ds.\nonumber
\end{equation}
It is easy to show that $\Upsilon(q,v)$ is a compact perturbation of the
identity such that for all $v\not\equiv 0$,
by definition of $\chi(q)$,
$\Upsilon(q,v)\neq0$, for all $q\in [2,p]$. Hence, by [\ref{De}, Theorem 8.10]
and the invariance of the degree under homo-topology, we have
\begin{equation}
\deg \left(I-T^p_\mu, B_r(0),0 \right)=\deg\left(I-T^2_\mu,
B_r(0),0\right)=\left\{
\begin{array}{l}
1, ~~\,\,\,\text{if}\,\,\mu\in
\left(0,\mu_1(p)\right),\\
-1,\,\,\text{if}\,\,\mu\in\left(\mu_1(p),\delta\right).
\end{array}
\right.\nonumber
\end{equation}\qed

\section{Unilateral Global bifurcation result}

\bigskip

\quad\, With a simple transformation $v=-u$, problem (\ref{Mg}) can be written as
\begin{equation}\label{B1}
\left\{
\begin{array}{l}
\left(\left(-v'\right)^N\right)'=\lambda^NNr^{N-1}\left(v^N+g(v)\right)\,\, \text{in}\,\, 0<r<1,\\
v'(0)=v(1)=0.
\end{array}
\right.
\end{equation}
Define the map $T_g:X\rightarrow X$ by
\begin{equation}
T_gv(r)=\int_r^1\left(\int_0^sN\tau^{N-1}\left((v(\tau))^N+g(v(\tau))\right)
\,d\tau\right)^{\frac{1}{N}}\,ds,\,\, 0\leq r\leq1.\nonumber
\end{equation}
It is not difficult to verify that $T_g$ is continuous and compact.
Clearly, problem (\ref{B1}) can be equivalently written as
\begin{equation}
v=\lambda T_g v.\nonumber
\end{equation}

\indent Now, we show that the existence and uniqueness theorem is valid for problem (\ref{B1}).
\\ \\
\textbf{Lemma 3.1.} \emph{If $(\lambda, v)$ is a solution of (\ref{B1})
and $v$ has a double zero, then $v \equiv 0$.}
\\ \\
\textbf{Proof.} Let $v$ be a solution of problem (\ref{B1}) and $r_*\in[0, 1]$ be a double zero.
We note that
\begin{equation}
v(r)=\lambda\int_r^{r_*}\left(\int_{r_*}^sN\tau^{N-1}\left((v(\tau))^N+g(v(\tau))\right)
\,d\tau\right)^{\frac{1}{N}}\,ds.\nonumber
\end{equation}
Firstly, we consider $r\in[0, r_*]$. Then we have
\begin{eqnarray}
\vert v(r)\vert&\leq&\lambda\left(\int_{r}^{r_*}N\tau^{N-1}\left\vert\left((v(\tau))^N
+g(v(\tau))\right)\right\vert\,d\tau\right)^{\frac{1}{N}},\nonumber
\end{eqnarray}
furthermore,
\begin{eqnarray}
\vert v(r)\vert^N&\leq&\lambda^N\int_{r}^{r_*}N\tau^{N-1}\left\vert\left((v(\tau))^N
+g(v(\tau))\right)\right\vert\,d\tau \nonumber\\
&\leq&\lambda^N\int_r^{r^*}
 N \tau^{N-1}\left\vert 1+\frac{g(v(\tau))}{(v(\tau))^N}\right\vert \vert v(\tau)\vert^N\,d\tau.\nonumber
\end{eqnarray}
According to the assumptions on $g$, for any $\varepsilon>0$, there exists a constant $\delta>0$ such that
\begin{equation}
\vert g(s)\vert\leq \varepsilon \vert s\vert^N\nonumber
\end{equation}
for any $\vert s\vert\in[0,\delta]$.
Hence, we have
\begin{equation}
\vert v(r)\vert^N\leq \lambda^N\int_r^{r^*}
N\left(1+\varepsilon+\max_{\vert s\vert\in\left[\delta,\Vert v\Vert\right]}
\left\vert\frac{g(s)}{s^N}\right\vert\right) \vert v(\tau)\vert^N\,d\tau.\nonumber
\end{equation}
By the modification of Gronwall-Bellman inequality [\ref{ILL}, Lemma 2.2], we get $v \equiv 0$ on $[0, r^*]$.
Similarly, using the Gronwall-Bellman inequality [\ref{Bre}, \ref{E}], we can get $v \equiv 0$ on $[r^*, 1]$
and the proof is complete.\qed\\

As Lions [\ref{L1}] showed, the first eigenvalue $\lambda_1$ is positive and
simple. Moreover, we also have the following result.\\
\\
\textbf{Lemma 3.2.} \emph{If $(\mu,\varphi)\in(0,+\infty)\times \left(C^2[0,1]
\setminus\{0\}\right)$ satisfies
\begin{equation}\label{B2}
\left\{
\begin{array}{l}
\left(\left(-v'\right)^N\right)'=\lambda^NNr^{N-1}v^N\,\, \text{in}\,\, 0<r<1,\\
v'(0)=v(1)=0
\end{array}
\right.
\end{equation}
 and $\mu\neq\lambda_1$,
then $\varphi$ must change sign.}\\
\\
\textbf{Proof.} By way of contradiction, we may suppose that $\varphi$ is not changing-sign. Without loss of
generality, we can assume that $\varphi\geq 0$ in
$(0,1)$. Lemma 3.1 follows that $\varphi>0$ in $(0,1)$. Theorem 1 of [\ref{L1}]
implies $\mu=\lambda_1$ and $\varphi=\theta \psi_1$
for some positive constant $\theta$, where $\psi_1$ is the positive eigenfunction
corresponding to $\lambda_1$ with $\Vert \psi_1\Vert=1$. We have a contradiction.\qed\\
\\
\indent Next, we show that $\lambda_1$ is also isolated.
\\ \\
\textbf{Lemma 3.3.} \emph{$\lambda_1$ is isolated; that is to say, $\lambda_1$ is
the unique eigenvalue in $(0,\delta)$ for some $\delta>\lambda_1$.}\\ \\
\textbf{Proof.} Theorem 1 of [\ref{L1}] has shown that $\lambda_1$ is left-isolated.
Assume by contradiction that there exists a sequence of eigenvalues
$\lambda_n\in(\lambda_1, \delta)$ which converge to $\lambda_1$. Let $v_n$ be the
corresponding eigenfunctions. Let $w_n:=v_n/\left\Vert v_n\right\Vert_{C^1[0,1]}$,
then $w_n$ should be the solutions of the problem
\begin{equation}
w_n=\lambda_n\int_r^1\left(\int_0^sN\tau^{N-1}w_n^N\,d\tau\right)^{\frac{1}{N}}\,ds.\nonumber
\end{equation}
Clearly, $w_n$ are bounded in $C^1[0,1]$ so there exists a subsequence, denoted again by $w_n$,
and $\psi\in X$ such that
$w_n\rightarrow\psi$ in $X$.
It follows that
\begin{equation}
\psi=\lambda_1\int_r^1\left(\int_0^sN\tau^{N-1}\psi^N\,d\tau\right)^{\frac{1}{N}}\,ds.\nonumber
\end{equation}
Then Theorem 1 of [\ref{L1}] follows that $\psi=\theta\psi_1$ for some positive constant $\theta$
in $(0,1)$. Thus $w_n\geq 0$ for $n$ large enough which contradicts $v_n$
changing-sign in $(0,1)$ which is implied by Lemma 3.2.\qed\\

Set
\begin{equation}
T_Nv:=\int_r^1\left(\int_0^sN\tau^{N-1}v^N\,d\tau\right)^{\frac{1}{N}}\,ds,\,\, 0\leq r\leq1.\nonumber
\end{equation}
Clearly, $I-T_N$ is a completely continuous vector field in
$X$. Thus, the Leray-Schauder degree
$\deg\left(I-T_N, B_r(0),0\right)$ is well defined for
arbitrary $r$-ball $B_r(0)$ and $\mu\in (0,\delta)$, where $\delta$ comes from Lemma 3.3.
\\ \\
\noindent\textbf{Lemma 3.4.} \emph{Let $\lambda$ be a constant with
$\lambda\in(0,\delta)$. Then for
arbitrary $r>0$,}
\begin{equation}
\deg \left(I-\lambda T_N, B_r(0),0 \right)=\left\{
\begin{array}{l}
1, \,\,\,\,\,\,\text{if}\,\,\lambda\in
\left(0,\lambda_1\right),\\
-1,\,\,\text{if}\,\,\lambda\in\left(\lambda_1,\delta\right).
\end{array}
\right.\nonumber
\end{equation}
\noindent\textbf{Proof.} Taking $p=N+1$ and $\mu=\lambda$ in $T_\mu^p$, we can see that
$\lambda_1=\mu_1(p)$. Furthermore, it is
no difficulty to verify that $T_\mu^p(v)=0$ for $\mu\in(0,\delta)$ implies that $v$ is not
changing-sign. It follows that $\lambda T_N=T_\mu^p$.
By Theorem 2.1, we can deduce this lemma.\qed
\\ \\
\textbf{Theorem 3.1.} \emph{$\left(\lambda_1,0\right)$ is a bifurcation
point of (\ref{B1}) and the associated bifurcation branch $\mathcal{C}$ in $\mathbb{R}\times X$
whose closure contains $\left(\lambda_1, 0\right)$ is either unbounded or contains a pair ($\overline{\lambda}, 0$)
where $\overline{\lambda}$ is an eigenvalue of (\ref{B2}) and $\overline{\lambda}\neq \lambda_1$.}
\\ \\
\textbf{Proof.} Suppose that $(\lambda_1, 0)$
is not a bifurcation point of problem (\ref{B1}). Then
there exist $\varepsilon > 0$, $\rho_0> 0$ such that for $\left\vert \lambda-\lambda_1\right\vert\leq\varepsilon$
and $0<\rho < \rho_0$ there is no
nontrivial solution of the equation
\begin{equation}
v-\lambda T_gv=0\nonumber
\end{equation}
with $\Vert v\Vert=\rho$. From the invariance of the degree under a compact
homotopy we obtain that
\begin{equation}\label{edc}
\text{deg}\left(I-\lambda T_gv,B_\rho(0),0\right)\equiv constant
\end{equation}
for $\lambda\in\left[\lambda_1-\varepsilon,\lambda_1+\varepsilon\right]$.

By taking $\varepsilon$ smaller if necessary, in view of Lemma 3.3, we can assume
that there is no eigenvalue
of (\ref{B2}) in $\left(\lambda_1,\lambda_1+\varepsilon\right]$. Fix $\lambda\in\left(\lambda_1,\lambda_1+\varepsilon\right]$.
We claim that the equation
\begin{equation}\label{es}
v-\lambda\int_r^1\left(\int_0^sN\tau^{N-1}\left(v^N+tg(v)\right)\,d\tau\right)^{\frac{1}{N}}\,ds=0
\end{equation}
has no solution $v$ with $\Vert v\Vert=\rho$ for every $t\in[0, 1]$ and $\rho$
sufficiently small.
Suppose on the contrary, let $\left\{v_n\right\}$ be the nontrivial solutions of (\ref{es})
with $\left\Vert v_n\right\Vert\rightarrow 0$
as $n\rightarrow+\infty$.

Let $w_n:=v_n/\left\Vert v_n\right\Vert$, then $w_n$ should be the solutions of the problem
\begin{equation}\label{B3}
w_n(t)=\lambda\int_r^1\left(\int_0^sN\tau^{N-1}\left(w_n^N+t\frac{g(v)}
{\left\Vert v_n\right\Vert^N}\right)\,d\tau\right)^{\frac{1}{N}}\,ds.
\end{equation}
Let
\begin{equation}
\widetilde{g}(v)=\max_{0\leq \vert s\vert\leq v}\vert g(s)\vert,\nonumber
\end{equation}
then $\widetilde{g}$ is nondecreasing with respect to $v$ and
\begin{equation}\label{eg0+}
\lim_{ v\rightarrow 0^+}\frac{\widetilde{g}(v)}{
v^{N}}=0.
\end{equation}
Further it follows from (\ref{eg0+}) that
\begin{equation}\label{egn0}
\frac{\vert g(v)\vert}{\Vert v\Vert^{N}} \leq \frac{\widetilde{ g}(v)}
{\Vert v\Vert^{N}}\leq \frac{
\widetilde{g}(\Vert v\Vert)}{\Vert v\Vert^{N}}\rightarrow0\,\,  \text{as}\,\, \Vert
v\Vert\rightarrow 0.
\end{equation}
By (\ref{B3}), (\ref{egn0}) and the compactness of $T_g$, we obtain that for
some convenient subsequence
$w_n\rightarrow w_0$ as $n\rightarrow+\infty$. Now $(\lambda,w_0)$ verifies
problem (\ref{B2}) and $\left\Vert w_0\right\Vert = 1$. This implies that $\lambda$ is an eigenvalue of (\ref{B2}).
This is a contradiction.

From the invariance of the degree under
homotopies and Lemma 3.4 we then obtain
\begin{equation}\label{edFk}
\deg\left(I-\lambda T_g(\cdot), B_r(0),0\right)=
\deg\left(I-\lambda T_N(\cdot), B_r(0),0\right)=-1.
\end{equation}
Similarly, for $\lambda\in \left[\lambda_1 - \varepsilon, \lambda_1\right)$ we find that
\begin{equation}\label{edFk1}
\deg\left(I-\lambda T_g(\cdot), B_r(0),0\right)=1.
\end{equation}
Relations (\ref{edFk}) and (\ref{edFk1}) contradict (\ref{edc}) and hence $\left(\lambda_1, 0\right)$ is a
bifurcation point of problem (\ref{B1}).

By standard arguments in global bifurcation theory (see [\ref{R}]), we can
show the existence of
a global branch of solutions of problem (\ref{B1}) emanating from
$\left(\lambda_1, 0\right)$. Our conclusion is proved.\qed\\

Next, we shall prove that the first choice of the alternative of Theorem 3.1 is
the only possibility. Let $P^+$ denote the set of functions
in $X$ which are positive in (0,1). Set $P^-=-P^+$ and $P =P^+\cup P^-$.
It is clear that $P^+$ and $P^-$ are disjoint and open in $X$.
Finally, let $K^{\pm}=\mathbb{R}\times P^{\pm}$ and $K=\mathbb{R}\times P$
under the product topology.
\\ \\
\textbf{Lemma 3.5.} \emph{The last alternative of Theorem 3.1 is
impossible if
$\mathcal{C}\subset \left(K\cup\{\left(\lambda_1,0\right)\}\right)$.}
\\ \\
\textbf{Proof.} Suppose on the contrary, if there exists $\left(\lambda_n,v_n\right)
\rightarrow\left(\overline{\lambda},0\right)$
when $n\rightarrow+\infty$ with $\left(\lambda_n,v_n\right)\in \mathcal{C}$,
$v_n \not\equiv 0$ and $\overline{\lambda}$
is another eigenvalue of (\ref{B2}).
Let $w_n :=v_n/\left\Vert v_n\right\Vert$, then $w_n$ should be the solutions of the problem
\begin{equation}\label{evs}
w_n=\lambda_n\int_r^1\left(\int_0^sN\tau^{N-1}\left(w_n^N+\frac{g(v)}
{\left\Vert v_n\right\Vert^N}\right)\,d\tau\right)^{\frac{1}{N}}\,ds.
\end{equation}
By an argument similar to that of Theorem 3.1, we obtain that for some convenient
subsequence $w_n\rightarrow w_0$ as $n\rightarrow+\infty$. It is easy to
see that $\left(\overline{\lambda}, w_0\right)$ verifies problem (\ref{B2})
and $\Vert w_0\Vert = 1$. Lemma 3.2 follows $w_0$ must change sign, and as a
consequence for some $n$ large enough, $w_n$ must change sign, and
this is a contradiction.\qed
\\ \\
\textbf{Remark 3.1.} Clearly, the proof of Lemma 3.5 also shows that
$\left(\lambda_1,0\right)$ is the unique bifurcation point
from $(\lambda,0)$ to the one-sign solutions of problem (\ref{B1}).
\\ \\
\textbf{Theorem 3.2.} \emph{There exists an unbounded
continuum $\mathcal{C}\subseteq K$ of solutions to problem (\ref{B1})
emanating from $\left(\lambda_1,0\right)$.}\\ \\
\textbf{Proof.}
Taking into account Theorem 3.1 and Lemma 3.5, we only need to prove that
$\mathcal{C}\subset \left(K\cup\{\left(\lambda_1,0\right)\}\right)$.
Suppose $\mathcal{C}\not\subset
\left(K\cup\left\{\left(\lambda_1,0\right)\right\}\right)$. Then there exists
$(\lambda, v)\in \left(\mathcal{C}\cap(\mathbb{R}\times \partial P)\right)$
such that $(\lambda, v) \neq
\left(\lambda_1, 0\right)$, $v\not\in P$, and $\left(\lambda_n,v_n\right)\rightarrow(\lambda, v)$
with $\left(\lambda_n,v_n\right)\in \left(\mathcal{C}
\cap(\mathbb{R}\times P)\right)$.
Since $v\in \partial P$, by Lemma 3.1, $v\equiv 0$. Let $w_n :=v_n/\left\Vert v_n\right\Vert$.
Using the proof similar to that of Lemma 3.5, we can show that
there exists $w\in X$ such that $(\lambda, w)$ satisfies (\ref{B2}) and $\Vert w\Vert = 1$,
that is to say, $\lambda$ is an eigenvalue of (\ref{B2}). Therefore, $\left(\lambda_n,v_n\right)\rightarrow
\left(\lambda, 0\right)$ with $\left(\lambda_n,v_n\right)\in \mathcal{C}\cap (\mathbb{R}\times P)$.
This contradicts Lemma 3.5.\qed
\\

Using an argument similar to one of [\ref{DM}, Theorem 3.2] with obvious changes, we
may obtain the following unilateral global bifurcation result.\\ \\
\textbf{Theorem 3.3.} \emph{There are two distinct unbounded sub-continua of solutions to problem (\ref{B1}),
$\mathcal{C}^+$ and $\mathcal{C}^-$ consisting of the bifurcation branch $\mathcal{C}$ and
\begin{equation}
\mathcal{C}^\nu\subset \left(K^\nu\cup\{\left(\lambda_1,0\right)\}\right),\nonumber
\end{equation}
where $\nu\in\{+,-\}$.}\\

\section{One-sign solutions}

\bigskip

\quad\, In this section, we shall investigate the existence and multiplicity of one-sign
solutions to problem (\ref{MO}). With a simple transformation $v=-u$, problem (\ref{MO})
can be written as
\begin{equation}\label{C1}
\left\{
\begin{array}{l}
\left(\left(-v'(r)\right)^N\right)'=\lambda^NNr^{N-1}f(v(r))\,\, \text{in}\,\, 0<r<1,\\
v'(0)=v(1)=0.
\end{array}
\right.
\end{equation}
Define the map $T_f:X\rightarrow X$ by
\begin{equation}
T_fv(r)=\int_r^1\left(\int_0^sN\tau^{N-1}f(v(r))\,d\tau\right)^{\frac{1}{N}}
\,ds,\,\, 0\leq r\leq1.\nonumber
\end{equation}
Similar to $T_g$, $T_f$ is continuous and compact. Clearly, problem (\ref{C1})
can be equivalently written as
\begin{equation}
v=\lambda T_f v.\nonumber
\end{equation}

\indent Let $f_0, f_\infty\in \mathbb{R}\setminus \mathbb{R}^-$ be such that
\begin{equation}
f_0^N=\lim_{ s\rightarrow0}\frac{f(s)}{s^N}\,\,\text{and}\,\,f_\infty^N
=\lim_{\vert s\vert\rightarrow+\infty}\frac{f(s)}{s^N}.\nonumber
\end{equation}
Through out this section, we always suppose that $f$ satisfies the following signum condition\\

(f1) \emph{$f\in C(\mathbb{R},\mathbb{R})$ with
$f(s)s^N>0$ for $s\neq0$.}\\

Clearly, (f1) implies $f(0)=0$. Hence, $v=0$ is always the solution of problem (\ref{MO}). Applying
Theorem 3.2, we shall establish the existence of one-sign
solutions of (\ref{MO}) as follows.\\ \\
\textbf{Theorem 4.1.} \emph{If $f_0\in(0,+\infty)$ and $f_\infty\in(0,+\infty)$,
then for any $\lambda\in\left(\lambda_1/f_\infty,\lambda_1/f_0\right)$
or $\lambda\in\left(\lambda_1/f_0,\lambda_1/f_\infty\right)$, (\ref{MO})
has two solutions $u^+$ and $u^-$ such
that $u^+$ is positive, strictly concave in $(0,1)$, and $u^-$ is negative,
strictly convex in $(0,1)$.}
\\ \\
\textbf{Proof.} It suffices to prove that (\ref{C1}) has two one-sign
solutions $v^+$ and $v^-$ such
that $v^+$ is positive, strictly concave in $(0,1)$, and $v^-$ is negative, strictly convex in $(0,1)$.

Let $\zeta\in C(\mathbb{R})$ be such that $f(s)=f_0^N s^N+\zeta(s)$
with $\lim_{s\rightarrow0}\zeta(s)/s^N=0.$
Applying Theorem 3.3 to (\ref{C1}), we have that
there are two distinct unbounded sub-continua, $\mathcal{C}^+$ and $\mathcal{C}^-$ consisting of the bifurcation branch $\mathcal{C}$ emanating
from $\left(\lambda_1/f_0, 0\right)$, such that
\begin{equation}
\mathcal{C}^\nu\subset \left(\left\{\left(\lambda_1,0\right)\right\}
\cup\left(\mathbb{R}\times P^\nu\right)\right).\nonumber
\end{equation}

To complete this theorem, it will be enough to show that $\mathcal{C}^\nu$ joins
$\left(\lambda_1/f_0, 0\right)$ to
$\left(\lambda_1/f_\infty, +\infty\right)$. Let
$\left(\mu_n, v_n\right) \in \mathcal{C}^\nu$ satisfy
$\mu_n+\left\Vert v_n\right\Vert\rightarrow+\infty.$
We note that $\mu_n >0$ for all $n \in \mathbb{N}$ since (0,0) is
the only solution of (\ref{C1}) for $\lambda= 0$ and
$\mathcal{C}^\nu\cap\left(\{0\}\times X\right)=\emptyset$.

We divide the rest proofs into two steps.

\emph{Step 1.} We show that there exists a constant $M$ such that $\mu_n
\in(0,M]$ for $n\in \mathbb{N}$ large enough.

On the contrary, we suppose that $\lim_{n\rightarrow +\infty}\mu_n=+\infty.$
On the other hand, we note that
\begin{equation}
\left(\left(-v_n'(r)\right)^N\right)'=\mu_n^N r^{N-1}
\widetilde{f}_n(r)v_n^N,\nonumber
\end{equation}
where
\begin{equation}
\widetilde{f}_n(r)=\left\{
\begin{array}{l}
\frac{f(v_n)}{v_n^N},\,\, \text{if}\,\,v_n\neq0,\\
f_0^N,\,\,\,\,\,\,\text{if}\,\,v_n=0.
\end{array}
\right.\nonumber
\end{equation}
The signum condition (f1) implies that there exists a positive
constant $\varrho$ such that $\widetilde{f}_n(r)\geq\varrho$ for
any $r\in[0,1]$.
By Lemma 3.2, we get $v_n$ must change sign in $(0,1)$ for $n$
large enough, and this contradicts the fact that $v_n\in \mathcal{C}^\nu$.

\emph{Step 2.} We show that $\mathcal{C}^\nu$ joins
$\left(\lambda_1/f_0, 0\right)$ to
$\left(\lambda_1/ f_\infty, +\infty\right)$.

It follows from \emph{Step 1} that $\left\Vert v_n\right\Vert\rightarrow+\infty.$
Let $\xi\in C(\mathbb{R})$ be such that $f(s)=f_\infty^N s^N+\xi(s).$
Then $\lim_{\vert s\vert\rightarrow+\infty}\xi(s)/s^N=0.$
Let $\widetilde{\xi}(v)=\max_{0\leq \vert s\vert\leq v}\vert \xi(s)\vert.$
Then $\widetilde{\xi}$ is nondecreasing and
\begin{equation}\label{eu0}
\lim_{v\rightarrow +\infty}\frac{\widetilde{\xi}(v)}{v^{N}}=0.
\end{equation}

We divide the equation
\begin{equation}
\left(\left(-v_n'\right)^N\right)'-\mu_n^Nf_\infty^N r^{N-1}v_n^N
=\mu_n^N r^{N-1}\xi(v_n)\nonumber
\end{equation}
by $\left\Vert v_n\right\Vert$ and set $\overline{v}_n = v_n/\left\Vert v_n\right\Vert$.
Since $\overline{v}_n$ are bounded in $X$,
after taking a subsequence if
necessary, we have that $\overline{v}_n \rightharpoonup \overline{v}$
for some $\overline{v} \in X$. Moreover, from
(\ref{eu0}) and the fact that $\widetilde{\xi}$ is nondecreasing,
we have that
\begin{equation}\label{C2}
\lim_{n\rightarrow+\infty}\frac{ \xi\left(v_n(r)\right)}{\left\Vert v_n\right\Vert^{N}}=0
\end{equation}
since
\begin{equation}
\frac{ \vert\xi\left(v_n(r)\right)\vert}{\left\Vert v_n\right\Vert^{N}}\leq\frac{ \widetilde{\xi}
(\left\vert v_n(r)\right\vert)}{\left\Vert v_n\right\Vert^{N}}
\leq\frac{ \widetilde{\xi}(\left\Vert v_n(r)\right\Vert)}{\left\Vert v_n\right\Vert^{N}}.\nonumber
\end{equation}

By the continuity and compactness of $T_f$, it follows that
\begin{equation}
\left(\left(-\overline{v}'\right)^N\right)'-\overline{\lambda}^Nf_\infty^N r^{N-1}
\overline{v}^N=0,\nonumber
\end{equation}
where
$\overline{\lambda}=\underset{n\rightarrow+\infty}\lim\lambda_n$, again
choosing a subsequence and relabeling it if necessary.

It is clear that $\Vert \overline{v}\Vert=1$ and $\overline{v}\in
\overline{\mathcal{C}^\nu}\subseteq
\mathcal{C}^\nu$ since $\mathcal{C}^\nu$ is closed in $\mathbb{R}\times X$.
Therefore, $\overline{\lambda}f_\infty=\lambda_1$, so
that $\overline{\lambda}=\lambda_1/ f_\infty.$
Therefore, $\mathcal{C}^\nu$ joins $\left(\lambda_1/
f_0, 0\right)$ to $\left(\lambda_1/f_\infty,
+\infty\right)$.\qed\\
\\
\textbf{Remark 4.1.} From the proof of Theorem 4.1, we can see
that if $f_0, f_\infty\in(0,+\infty)$ then there exist $\lambda_2^\nu>0$ and $\lambda_3^\nu>0$ such that
(\ref{MO}) has at least a strictly convex solution or a strictly
concave solution for all $\lambda\in\left(\lambda_2^\nu,\lambda_3^\nu\right)$
and has no nontrivial convex or concave solution for all
$\lambda\in\left(0,\lambda_2^\nu\right)\cup\left(\lambda_3^\nu,+\infty\right)$.
\\ \\
\textbf{Proof.} It is sufficient to show that there exist $\lambda_2^\nu>0$ such that
(\ref{MO}) has no nontrivial convex or concave solution for all
$\lambda\in\left(0,\lambda_2^\nu\right)$. Suppose on the contrary that there exists a sequence $\left\{\mu_n, v_n\right\}\in\mathcal{C}^\nu$ such that
$\underset{n\rightarrow +\infty}{\lim}\mu_n=0$ and $v_n\not\equiv 0$. $f_0, f_\infty\in(0,+\infty)$ implies that there exists a positive constant $M$ such that
\begin{equation}
\left\vert\frac{f(s)}{s^N}\right\vert\leq M \,\,\text{for any}\,\, s\neq 0.\nonumber
\end{equation}
Let $w_n =v_n/\left\Vert v_n\right\Vert$. Obviously, one has
\begin{eqnarray}
1=\left\Vert w_n\right\Vert=\left\Vert\mu_n\int_r^1\left(\int_0^sN\tau^{N-1}\left(\frac{f\left(v_n\right)}
{\left\Vert v_n\right\Vert^N}\right)\,d\tau\right)^{\frac{1}{N}}\,ds\right\Vert
\leq M^{1/N}\mu_n\rightarrow 0. \nonumber
\end{eqnarray}
This is a contradiction.\qed\\

From the proof of Theorem 4.1 and  Remark 4.1, we can deduce the following two corollaries.
\\
\\
\textbf{Corollary 4.1.} \emph{Assume that there exists a positive constant $\rho>0$ such that
\begin{equation}
\frac{f(s)}{s^N}\geq \rho\nonumber
\end{equation}
for any $s\neq 0$.
Then there exist $\zeta_*^+>0$ and $\zeta_*^-<0$ such that problem (\ref{MO}) has no one-sign solution for
any $\lambda\in\left(-\infty,\zeta_*^-\right)\cup\left(\zeta_*^+,+\infty\right)$.}
\\ \\
\textbf{Corollary 4.2.} \emph{Assume that there exists a positive constant $\varrho>0$ such that
\begin{equation}
\left\vert\frac{f(s)}{s^N}\right\vert\leq\varrho\nonumber
\end{equation}
for any $s\neq 0$.
Then there exist $\eta_*^+>0$ and $\eta_*^-<0$ such that problem (\ref{MO}) has no one-sign solution for
any $\lambda\in\left(0,\eta_*^-\right)\cup\left(0,\eta_*^+\right)$.}
\\ \\
\textbf{Theorem 4.2.} \emph{If $f_0\in(0,+\infty)$ and $f_\infty=0$,
then for any $\lambda\in\left(\lambda_1/f_0,+\infty\right)$, (\ref{MO})
has two solutions $u^+$ and $u^-$ such
that $u^+$ is positive, strictly concave in $(0,1)$, and $u^-$ is negative,
strictly convex in $(0,1)$.}
\\ \\
\textbf{Proof.} In view of Theorem 4.1, we only need to show that
$\mathcal{C}^\nu$ joins
$\left(\lambda_1/f_0, 0\right)$ to
$\left(+\infty, +\infty\right)$.
Suppose on the contrary that there exists $\mu_M$ be a blow up point
(see Definition 1.1 of [\ref{S}]) and $\mu_M<+\infty$.
Then there exists a sequence $\left\{\mu_n, v_n\right\}$ such that
$\underset{n\rightarrow +\infty}{\lim}\mu_n=\mu_{M}$ and
$\underset{n\rightarrow +\infty}{\lim}\left\Vert v_n\right\Vert=+\infty$ as
$n\rightarrow+\infty$. Let $w_n =v_n/\left\Vert v_n\right\Vert$ and $w_n$
should be the solutions of the problem
\begin{equation}
w_n=\mu_n\int_r^1\left(\int_0^sN\tau^{N-1}\left(\frac{f\left(v_n\right)}
{\left\Vert v_n\right\Vert^N}\right)\,d\tau\right)^{\frac{1}{N}}\,ds.\nonumber
\end{equation}
Similar to (\ref{C2}), we can show
\begin{equation}
\lim_{n\rightarrow+\infty}\frac{ f\left(v_n(r)\right)}{\left\Vert v_n\right\Vert^{N}}=0.\nonumber
\end{equation}
By the compactness of $T_f$, we obtain that for some convenient subsequence
$w_n\rightarrow w_0$ as $n\rightarrow+\infty$. Letting $n\rightarrow+\infty$,
we obtain that
$w_0\equiv0$. This contradicts $\left\Vert w_0\right\Vert=1$. \qed\\
\\
\textbf{Remark 4.2.} Under the assumptions of Theorem 4.2, in view of Corollary 4.2, we can see that there exists $\lambda_4^\nu>0$ such that
problem (\ref{MO}) has at least a strictly convex solution or a strictly concave
solution for all $\lambda\in\left(\lambda_4^\nu,+\infty\right)$
and has no nontrivial convex or concave solution for all $\lambda\in\left(0,\lambda_4^\nu\right)$.
\\ \\
\textbf{Theorem 4.3.} \emph{If $f_0\in(0,+\infty)$ and $f_\infty=\infty$,
then for any $\lambda\in\left(0,\lambda_1/f_0\right)$, (\ref{MO}) has two
solutions $u^+$ and $u^-$ such
that $u^+$ is positive, strictly concave in $(0,1)$, and $u^-$ is negative,
strictly convex in $(0,1)$.}
\\ \\
\textbf{Proof.} Considering of the proof of Theorem 4.1, we only need to show that $\mathcal{C}^\nu$ joins
$\left(\lambda_1/f_0, 0\right)$ to
$\left(0, +\infty\right)$. Clearly, $f_\infty=+\infty$ implies that $f(s)\geq M^Ns^N$
for some positive constant $M$ and $\vert s\vert$ large enough.

To complete the proof, it suffices to show that the
unique blow up point of $\mathcal{C}^\nu$ is $\lambda=0$. Suppose on the contrary
that there exists $0<\widehat{\lambda}$ is a
blow up point of $\mathcal{C}^\nu$. Then there exists a sequence $\left\{\lambda_n, v_n\right\}$
such that $\underset{n\rightarrow +\infty}{\lim}\lambda_n=\widehat{\lambda}$ and
$\underset{n\rightarrow +\infty}{\lim}\left\Vert v_n\right\Vert=+\infty$. Let $w_n =v_n/\left\Vert v_n\right\Vert$. Clearly, one has
\begin{equation}
w_n=\lambda_n\int_r^1\left(\int_0^sN\tau^{N-1}\left(\frac{f\left(v_n\right)}
{v_n^N}\frac{v_n^N}{\left\Vert v_n\right\Vert^N}\right)\,d\tau\right)^{\frac{1}{N}}\,ds.\nonumber
\end{equation}
Take $M=64/\widehat{\lambda}+1$. For $r\in\left[1/4,3/4\right]$, by virtue of Lemma 2.2
of [\ref{HW}], we have
\begin{eqnarray}\label{C3}
\left\vert w_n\right\vert&\geq& M\lambda_n\int_r^1
\left(\int_0^sN\tau^{N-1}\left\vert w_n\right\vert^N\,d\tau\right)^{\frac{1}{N}}\,ds\nonumber\\
&\geq& M\left\Vert w_n\right\Vert\lambda_n\int_r^1
\left(\int_0^sN\tau^{N-1}(1-\tau)^N\,d\tau\right)^{\frac{1}{N}}\,ds\nonumber\\
&\geq& M\left\Vert w_n\right\Vert\lambda_n(1-r)
\left(\int_0^rN\tau^{N-1}(1-\tau)^N\,d\tau\right)^{\frac{1}{N}}\nonumber\\
&\geq& M\left\Vert w_n\right\Vert\lambda_n(1-r)^{2}
\left(\int_0^rN\tau^{N-1}\,d\tau\right)^{\frac{1}{N}}\nonumber\\
&\geq& M\left\Vert w_n\right\Vert\lambda_nr(1-r)^{2}\nonumber\\
&\geq& \frac{M\left\Vert w_n\right\Vert\lambda_n}{64}.
\end{eqnarray}
It is obvious that (\ref{C3}) follows $M\lambda_n\leq 64$. Thus, we
get $M\leq 64/\widehat{\lambda}$. While, this is impossible because
of $M=64/\widehat{\lambda}+1$.\qed
\\ \\
\textbf{Remark 4.3.} Clearly, Theorem 4.3 and Corollary 4.1 imply that
if $f_0\in(0,+\infty)$ and $f_\infty=+\infty$ then there exists $\lambda_5^\nu>0$ such that
(\ref{MO}) has at least a strictly convex solution or a strictly
concave solution for all $\lambda\in\left(0,\lambda_5^\nu\right)$
and has no nontrivial convex or concave solution for all $\lambda\in\left(\lambda_5^\nu,+\infty\right)$.
\\ \\
\textbf{Theorem 4.4.} \emph{If $f_0=0$ and $f_\infty\in(0,+\infty)$,
then for any $\lambda\in\left(\lambda_1/f_\infty,+\infty\right)$, (\ref{MO})
has two solutions $u^+$ and $u^-$ such
that $u^+$ is positive, strictly concave in $(0,1)$, and $u^-$ is negative,
strictly convex in $(0,1)$.}
\\ \\
\textbf{Proof.} If $(\lambda,v)$ is any solution of (\ref{C1}) with
$\Vert v\Vert\not\equiv 0$, dividing (\ref{C1}) by $\Vert v\Vert^{2N}$ and
setting $w=v/\Vert v\Vert^2$ yields
\begin{equation}
\left\{
\begin{array}{l}
\left(\left(-w'(r)\right)^N\right)'=\lambda^NNr^{N-1}\left(\frac{f(v)}
{\Vert v\Vert^{2N}}\right)\,\, \text{in}\,\, 0<r<1,\\
w'(0)=w(1)=0.
\end{array}
\right.\label{C4}
\end{equation}
Define
\begin{equation}
\widetilde{f}(w)=\left\{
\begin{array}{l}
\Vert w\Vert^{2N}f\left(\frac{w}{\Vert w\Vert^2}\right),\,\,\text{if}\, w\neq 0,\\
0,~~~~~~~~~~~~~~~~~~~\,\, \text{if}\,\, w=0.
\end{array}
\right.\nonumber
\end{equation}

Clearly, (\ref{C4}) is equivalent to:
\begin{equation}\label{C5}
\left\{
\begin{array}{l}
\left(\left(-w'(r)\right)^N\right)'=\lambda^Nr^{N-1}\widetilde{f}(w)\,\, \text{in}\,\, 0<r<1,\\
w'(0)=w(1)=0.
\end{array}
\right.
\end{equation}
It is obvious that $(\lambda,0)$ is always the solution of (\ref{C5}).
By simple computation, we can show that $\widetilde{f}_0=f_\infty$
and $\widetilde{f}_\infty=f_0$.

Now applying Theorem 4.2 and the inversion $w\rightarrow w/\Vert w\Vert^2=v$,
we can achieve our conclusion.\qed\\
\\
\textbf{Remark 4.4.} Under the assumptions of Theorem 4.4, we note there exists
$\lambda_6^\nu>0$ such that
(\ref{MO}) has at least a strictly convex solution or a strictly concave solution
for all $\lambda\in\left(\lambda_6^\nu,+\infty\right)$
and has no nontrivial convex or concave solution for all $\lambda\in\left(0,\lambda_6^\nu\right)$.\\

\indent Next, we shall need the following topological lemma:
\\ \\
\textbf{Lemma 4.1 (see [\ref{MA}].} \emph{Let $X$ be a Banach space and let $C_n$ be a family of closed connected subsets of $X$. Assume that:}
\\

(i) \emph{there exist $z_n\in C_n$, $n=1,2,\ldots$, and $z^*\in X$, such that $z_n\rightarrow z^*$;}

(ii) \emph{$r_n=\sup \left\{\Vert x\Vert\big| x\in C_n\right\}=+\infty$;}

(iii) \emph{for every $R>0$, $\left(\cup_{n=1}^{+\infty} C_n\right)\cap B_R$ is a relatively compact set of $X$, where}
\begin{equation}
B_R=\{x\in X|\Vert x\Vert\leq R\}.\nonumber
\end{equation}

\noindent \emph{Then there exists an unbounded component $\mathfrak{C}$ in $\mathfrak{D} =\limsup_{n\rightarrow +\infty}C_n$ and $z\in \mathfrak{C}$.}
\\ \\
\textbf{Theorem 4.5.} \emph{If $f_0=0$ and $f_\infty=0$,
then there exists $\lambda_*^+>0$ such that for any $\lambda\in\left(\lambda_*^+,+\infty\right)$,
(\ref{MO}) has two solutions $u_1^+$ and $u_2^+$ such
that they are positive, strictly concave in $(0,1)$. Similarly, there exists $\lambda_*^->0$
such that for any $\lambda\in\left(\lambda_*^-,+\infty\right)$, (\ref{MO}) has two
solutions $u_1^-$ and $u_2^-$ such that they are negative, strictly convex in $(0,1)$.}
\\ \\
\textbf{Proof.} Define
\begin{equation}
f^n(s)=\left\{
\begin{array}{l}
\frac{1}{n^N}s^N,\,\,\quad\quad\quad\quad\quad\quad\quad\quad\quad\quad\quad\,\,\,\,\,
\quad\quad\quad\,\,\quad s\in\left[-\frac{1}{n},\frac{1}{n}\right],\\
\left(f\left(\frac{2}{n}\right)-\frac{1}{n^{2N}}\right)ns+\frac{2}{n^{2N}}
-f\left(\frac{2}{n}\right),\quad\quad\quad\quad\,\,\,\,\,\,s\in\left(\frac{1}{n},\frac{2}{n}\right),\\
-\left(f\left(-\frac{2}{n}\right)-\frac{(-1)^N}{n^{2N}}\right)ns+\frac{2(-1)^N}{n^{2N}}
-f\left(-\frac{2}{n}\right),\,\,s\in\left(-\frac{2}{n},-\frac{1}{n}\right),\\
f(s),\quad\quad\quad\quad\quad\quad\quad\quad\quad\quad\quad\quad\quad\quad\quad\,\,
\quad\,\,\,\,s\in\left(-\infty,-\frac{2}{n}\right]\cup\left[\frac{2}{n},+\infty\right).
\end{array}
\right.\nonumber
\end{equation}
Now, consider the following problem
\begin{equation}
\left\{
\begin{array}{l}
\left(\left(u'(r)\right)^N\right)'=\lambda^NNr^{N-1} f^n(-u(r))\,\, \text{in}\,\, 0<r<1,\\
u'(0)=u(1)=0.
\end{array}
\right.\nonumber
\end{equation}
Clearly, we can see that $\lim_{n\rightarrow+\infty}f^n(s)=f(s)$, $f_0^n=1/n$ and $f_\infty^n=f_\infty=0$.
Theorem 4.2 implies
that there exists a sequence unbounded continua $\mathcal{C}_n^\nu$ emanating from $\left(n\lambda_1, 0\right)$
and joining to $\left(+\infty, +\infty\right)$.

Taking $z_n^1=(n\lambda_1, 0)$ and $z_n^2=(+\infty, +\infty)$, we have $z_n^1, z_n^2\in \mathcal{C}_n^\nu$ and $z_n^1\rightarrow(+\infty,0)$, $z_n^2\rightarrow (+\infty,+\infty)$. The compactness of $T_f$ implies that $\left(\cup_{n=1}^{+\infty} \mathcal{C}_n^\nu\right)\cap B_R$ is pre-compact. So Lemma 4.1 implies that there exists an unbounded component $\mathcal{C}^\nu$ of $\limsup_{n\rightarrow +\infty}\mathcal{C}_n^\nu$ such that
$(+\infty,0)\in\mathcal{C}^\nu$ and $\left(+\infty, +\infty\right)\in \mathcal{C}^\nu$. By an argument similar to that of Theorem 4.2, we can show that $\mathcal{C}^\nu\cap([0,+\infty)\times\{0\})=\emptyset$.\qed\\
\\
\textbf{Remark 4.5.} From Theorem 4.5 and Corollary 4.2, we also can see that there exists $\lambda_7^\nu>0$ such that
(\ref{MO}) has at least two strictly convex solution or two strictly concave solution for all
$\lambda\in\left[\lambda_7^\nu,\lambda_*^\nu\right]$
and has no nontrivial convex or concave solution for all $\lambda\in\left(0,\lambda_7^\nu\right)$.
\\ \\
\textbf{Theorem 4.6.} \emph{If $f_0=0$ and $f_\infty=\infty$,
then for any $\lambda\in\left(0,+\infty\right)$, (\ref{MO}) has two solutions $u^+$ and $u^-$ such
that $u^+$ is positive, strictly concave in $(0,1)$, and $u^-$ is negative, strictly convex in $(0,1)$.}
\\ \\
\textbf{Proof.} Using an argument similar to that of Theorem 4.5, in view of the conclusion of
Theorem 4.3, we can easily get the results of this theorem.\qed
\\ \\
\textbf{Theorem 4.7.} \emph{If $f_0=\infty$ and $f_\infty=0$,
then for any $\lambda\in\left(0,+\infty\right)$, (\ref{MO}) has two solutions $u^+$ and $u^-$ such
that $u^+$ is positive, strictly concave in $(0,1)$, and $u^-$ is negative, strictly convex in $(0,1)$.}
\\ \\
\textbf{Proof.} By an argument similar to that of Theorem 4.4 and the conclusions of Theorem 4.6, we can prove it.\qed
\\ \\
\textbf{Theorem 4.8.} \emph{If $f_0=\infty$ and $f_\infty\in(0,+\infty)$,
then for any $\lambda\in\left(0,\lambda_1/f_\infty\right)$, (\ref{MO}) has
two solutions $u^+$ and $u^-$ such
that $u^+$ is positive, strictly concave in $(0,1)$, and $u^-$ is negative, strictly convex in $(0,1)$.}
\\ \\
\textbf{Proof.} By an argument similar to that of Theorem 4.4 and the conclusion of Theorem 4.3, we can obtain it.\qed
\\ \\
\textbf{Remark 4.6.} Similarly to Remark 4.3, there exists $\lambda_8^\nu>0$ such that
(\ref{MO}) has at least a strictly convex solution or a strictly concave solution
for all $\lambda\in\left(0,\lambda_8^\nu\right)$
and has no nontrivial convex or concave solution for all $\lambda\in\left(\lambda_8^\nu,+\infty\right)$.
\\ \\
\textbf{Theorem 4.9.} \emph{If $f_0=\infty$ and $f_\infty=\infty$,
then there exists $\mathcal{\lambda}^+>0$ such that for any
$\lambda\in\left(0,\lambda^+\right)$, (\ref{MO}) has two solutions $u_1^+$ and $u_2^+$ such
that they are positive, strictly concave in $(0,1)$. Similarly, there exists $\lambda^->0$
such that for any $\lambda\in\left(0,\lambda^-\right)$, (\ref{MO}) has two solutions $u_1^-$
and $u_2^-$ such that they are negative, strictly convex in $(0,1)$.}\\
\\
\textbf{Proof.} Define
\begin{equation}
f^n(s)=\left\{
\begin{array}{l}
n^Ns^N,\,\,\quad\quad\quad\quad\quad\quad\quad\quad\quad\quad\quad\,\,\,
\quad\quad\quad\quad\,\,\quad s\in\left[-\frac{1}{n},\frac{1}{n}\right],\\
\left(f\left(\frac{2}{n}\right)-1\right)ns+2-f\left(\frac{2}{n}\right),\quad\quad\quad\quad\quad\quad\quad\,\,\,\,\,s\in
\left(\frac{1}{n},\frac{2}{n}\right),\\
-\left(f\left(-\frac{2}{n}\right)-(-1)^N\right)ns+2(-1)^N-f\left(-\frac{2}{n}\right),\,\,s\in
\left(-\frac{2}{n},-\frac{1}{n}\right),\\
f(s),\quad\quad\quad\quad\quad\quad\quad\quad\quad\quad\quad\quad\,
\quad\quad\quad\quad\quad\,\,\,\,s\in\left(-\infty,-\frac{2}{n}\right]\cup\left[\frac{2}{n},+\infty\right).
\end{array}
\right.\nonumber
\end{equation}
By the conclusions of Theorem 4.3 and an argument similar to that of Theorem 4.5, we can prove there exists an unbounded component $\mathcal{C}^\nu$ of solutions to problem (\ref{MO}) such that $(0,0)\in\mathcal{C}^\nu$ and $\left(0, +\infty\right)\in \mathcal{C}^\nu$. By an argument similar to that of Theorem 4.3, we can show that $\mathcal{C}^\nu\cap((0,+\infty)\times\{0\})=\emptyset$. By arguments similar to those of Theorem 4.3 and 4.5, we can show that there exists $\mu_*^\nu>0$ such that $\mathcal{C}^\nu\cap\left(\left(\mu_*^\nu,+\infty\right)\times X\right)=\emptyset$.\qed
\\ \\
\textbf{Remark 4.7.} By Theorem 4.9 and Corollary 4.1, we can see that there exists $\lambda_9^\nu>0$ such that
(\ref{MO}) has at least a strictly convex solution or a
strictly concave solution for all $\lambda\in\left[\lambda^\nu,\lambda_9^\nu\right]$
and has no nontrivial convex or concave solution for all
$\lambda\in\left(\lambda_9^\nu,+\infty\right)$. \\
\\
\textbf{Remark 4.8.} Clearly, the conclusions of Theorem 1.1 of [\ref{W}] and Theorem 5.1 of [\ref{HW}]
are the corollaries of Theorem 4.1--4.9.\\
\\
\textbf{Remark 4.9.} Let $f(s)=e^s$. It can be easily verified that $f_0=\infty$ and $f_\infty=\infty$.
This fact with Remark 4.7 implies that there is no solution of problem (\ref{MO}) with $\lambda$ large enough,
and for sufficiently small $\lambda$ there are two strictly convex solutions. Set $\mu:=\lambda^{1/2}$.
Through a scaling, we can show that problem (\ref{MO}) is equivalent to
\begin{equation}\label{MAh}
\left\{
\begin{array}{l}
\det\left(D^2u\right)=e^{-u}\,\, \text{in}\,\, B_\mu(0),\\
u=0~~~~~~~~~~~~~~\,\text{on}\,\, \partial B_\mu(0),
\end{array}
\right.
\end{equation}
where $B_\mu(0)$ denotes the set of $\{x\in \mathbb{R}^N\big|\vert x\vert\leq \mu\}$.
Hence there is no solution of problem (\ref{MAh}) with $\mu$ large enough,
and for sufficiently small $\mu$ there are two strictly convex solutions. Obviously,
this result improve the corresponding one of [\ref{ZW}, Theorem 3.1]. So Theorem 3.1 of [\ref{ZW}] is
our corollary of Theorem 4.9.\\
\\
\textbf{Remark 4.10.} Obviously, the results of Theorem 4.1--4.9 are also valid on
$B_R(0)$ for any $R>0$.

\section{Exact multiplicity of one-sign solutions}

\bigskip

\quad\, In this section, under some more strict assumptions of $f$, we shall
show that the unbounded continuum which are obtained in Section 4
may be smooth curves. We just show the case of $f_0\in(0,+\infty)$ and $f_\infty=0$. Other cases are similar.

Firstly, we study the local structure of the bifurcation branch $\mathcal{C}$
near $\left(\lambda_1,0\right)$, which is obtained in Theorem 3.1.
Let $\mathbb{E}=\mathbb{R}\times X$, $\Phi(\lambda,v):=v-\lambda T_g(v)$ and
\begin{equation}
\mathcal{S}:=\overline{\left\{(\lambda,v)\in \mathbb{E}\big|\Phi(\lambda,v)=0,
v\neq0\right\}}^{\mathbb{E}}.\nonumber
\end{equation}
In order to formulate and prove main results of this section, it is convenient
to introduce L\'{o}pez-G\'{o}mez's notations [\ref{LG}].
Given any $\lambda\in \mathbb{R}$ and $0 < s < +\infty$, we consider an
open neighborhood of $\left(\lambda_1, 0\right)$ in $\mathbb{E}$ defined by
\begin{equation}
\mathbb{B}_s(\lambda_1,0):=\left\{(\lambda,v)\in \mathbb{E}\big|\Vert v\Vert
+\left\vert\lambda-\lambda_1\right\vert<s\right\}.\nonumber
\end{equation}
Let $X_0$ be a closed subspace of $X$ such that
\begin{equation}
X=\text{span}\left\{\psi_1\right\}\oplus X_0.\nonumber
\end{equation}
According to the Hahn-Banach theorem,
there exists a linear functional $l\in X^*$, here $X^*$ denotes the dual
space of $X$, such that
\begin{equation}
l\left(\psi_1\right)=1\,\, \text{and}\,\, X_0=\{v\in X\big|l(v)=0\}.\nonumber
\end{equation}
Finally, for any $0<\varepsilon<+\infty$ and $0 < \eta < 1$,
we define
\begin{equation}
K_{\varepsilon,\eta}:=\left\{(\lambda,v)\in \mathbb{E}\big| \left\vert\lambda-\lambda_1
\right\vert<\varepsilon, \vert l(v)\vert>\eta\Vert v\Vert\right\}.\nonumber
\end{equation}
Since
\begin{equation}
u\mapsto \vert l(u)\vert-\Vert u\Vert\nonumber
\end{equation}
is continuous, $K_{\varepsilon,\eta}$ is an open subset of $\mathbb{E}$
consisting of two disjoint components $K_{\varepsilon,\eta}^+$ and
$K_{\varepsilon,\eta}^-$, where
\begin{equation}
K_{\varepsilon,\eta}^+:=\left\{(\lambda,v)\in \mathbb{E}\big|  \left\vert\lambda-
\lambda_1\right\vert<\varepsilon, l(v)>\eta\Vert v\Vert\right\},\nonumber
\end{equation}
\begin{equation}
K_{\varepsilon,\eta}^-:=\left\{(\lambda,v)\in \mathbb{E}\big|  \left\vert\lambda-
\lambda_1\right\vert<\varepsilon, l(v)<-\eta\Vert v\Vert\right\}.\nonumber
\end{equation}

Applying an argument similar to prove [\ref{LG}, Lemma
6.4.1] with obvious changes, we may obtain the following result, which localizes
the possible solutions of (\ref{Mg})
bifurcating from $\left(\lambda_1,0\right)$.\\ \\
\textbf{Lemma 5.1.} \emph{For every $\eta\in(0, 1)$ there exists a number
$\delta_0>0$ such that
for each $0<\delta<\delta_0$,
\begin{equation}
\left(\left(\mathcal{S}\setminus\left\{\left(\lambda_1,0\right)\right\}\right)\cap \mathbb{B}_\delta
\left(\lambda_1,0\right)\right)\subset K_{\varepsilon,\eta}.\nonumber
\end{equation}
Moreover, for each
\begin{equation}
(\lambda,v)\in\left(\mathcal{S}\setminus\left\{\left(\lambda_1,0\right)\right\}\right)\cap
\left(\mathbb{B}_\delta\left(\lambda_1,0\right)\right),\nonumber
\end{equation}
there are $s\in\mathbb{R}$ and unique $y\in X_0$ such that
\begin{equation}
v=s\psi_1+y\,\, \text{and} \,\, \vert s\vert>\eta\Vert v\Vert.\nonumber
\end{equation}
Furthermore, for these solutions $(\lambda,v)$,
\begin{equation}
\lambda=\lambda_1+o(1)\,\, \text{and}\,\, y=o(s)\nonumber
\end{equation}
as $s\rightarrow 0$.}\\

\indent Moreover, the next lemma shows that the component $\mathcal{C}$ of
$S$ emanating from $\left(\lambda_1,0\right)$
consists of two sub-continua meeting each other at $\left(\lambda_1,0\right)$.\\
\\
\textbf{Lemma 5.2.} \emph{Let $\mathcal{C}$ denote the component of
$\mathcal{S}$ emanating from $\left(\lambda_1,0\right)$.
Then, $\mathcal{C}$ possesses a sub-continuum in each of the cones
\begin{equation}
K_{\varepsilon,\eta}^+\cup\left\{\left(\lambda_1,0\right)\right\}\,\,\text{and}\,\,
K_{\varepsilon,\eta}^-\cup\left\{\left(\lambda_1,0\right)\right\}\nonumber
\end{equation}
each of which meets $\left(\lambda_1,0\right)$ and $\partial B_\delta\left(\lambda_1,0\right)$ for all $\delta>0$
that are sufficiently small.}\\
\\
\textbf{Proof.} It is easy to show that $v=l(v)\psi_1+y$. We define
\begin{equation}\label{e1.1}
\widehat{g}(v)=\left\{
\begin{array}{l}
g(v)\,\,\quad\quad\quad\quad\quad\quad\quad\,\, \text{if\,\,}l(v)\leq -\eta\Vert v\Vert;\\
\frac{-l(v)}{\eta \Vert v\Vert}g\left(-\eta\Vert v\Vert\psi_1+y\right)\,\,
\text{if\,\,}-\eta\Vert v\Vert<l(v)\leq0;\\
-g(-v)\, \quad\quad\quad\quad\quad\quad\text{if\,\,}l(v)>0
\end{array}
\right.\nonumber
\end{equation}
and
\begin{equation}
\widehat{\Phi}(\lambda,v)=v-\lambda\int_r^1\left(\int_0^sN\tau^{N-1}
\left(v^N+\widehat{g}(v)\right)\,d\tau\right)^{\frac{1}{N}}\,ds.\nonumber
\end{equation}
Clearly, the mapping $\widehat{\Phi}(\lambda,v)$ is odd with respect to $v$.
Since the rest proof is similar to [\ref{LG}, Proposition 6.4.2], we omit it here.\qed\\
\\
\textbf{Remark 5.1.} From Lemma 5.1 and 5.2, we can see that $\mathcal{C}$ near
$\left(\lambda_1,0\right)$ is given by a curve
$(\lambda(s),v(s))=\left(\lambda_1+o(1),s\psi_1+o(s)\right)$ for $s$ near $0$.
Moreover, we can distinguish between two
portions of this curve by $s\geq 0$ and $s\leq 0$.
\\

\indent The primary result in this section is the following result.
\\ \\
\textbf{Theorem 5.1.} \emph{Let $f\in C^1(\mathbb{R})$ satisfies the
assumptions of Theorem 4.2. Suppose $f'(s)<Nf(s)/s$ for any $s>0$
and $f'(s)>Nf(s)/s$ for any $s<0$.
Then for any $\lambda\in\left(\lambda_1/f_0,+\infty\right)$, (\ref{MO})
has exactly two solutions $u^+$ and $u^-$ such
that $u^+$ is positive, strictly concave in $(0,1)$, and $u^-$ is negative,
strictly convex in $(0,1)$.}
\\ \\
\textbf{Remark 5.2.} Clearly, the assumption $f'(s)<Nf(s)/s$ for $s>0$ is equivalent
to $f(s)/s^N$ is decreasing for $s>0$. However, if $N$ is even (or odd) then
$f'(s)>Nf(s)/s$ for $s<0$ is equivalent to $f(s)/s^N$ is increasing (or decreasing) for $s<0$.
\\ \\
\indent We use the stability properties to prove Theorem 5.1.
Let
\begin{equation}
Y:=\left\{v\in C^2(0,1)\big|v'(0)=v(1)=0\right\}.\nonumber
\end{equation}
For any $\phi\in Y$ and one-sign solution $u$ of (\ref{MO}),
by some simple computations, we can show
that the linearized equation of (\ref{MO}) about $u$ at the direction $\phi$ is
\begin{equation}\label{EM1}
\left\{
\begin{array}{l}
\left(-\phi'\left(-v'\right)^{N-1}\right)'-\lambda^Nr^{N-1}f'(v)\phi=\frac{\mu}{N} \phi\,\, \text{in}
\,\,(0,1),\\
\phi'(0)=\phi(1)=0,
\end{array}
\right.
\end{equation}
where $v=-u$. Hence, the linear stability of a solution $u$ of (\ref{MO})
can be determined by the linearized eigenvalue problem (\ref{EM1}).
A solution $u$ of (\ref{MO}) is stable if all eigenvalues of (\ref{EM1})
are positive, otherwise it is unstable.
We define the \emph{Morse index} $M(u)$ of a solution $u$ to (\ref{MO})
to be the number of negative eigenvalues of (\ref{EM1}).
A solution $u$ of (\ref{MO}) is degenerate if $0$ is an eigenvalue of
(\ref{EM1}), otherwise it is non-degenerate. \\

\indent The following lemma is our main stability result for the negative steady state solution.\\ \\
\textbf{Lemma 5.3.} \emph{Suppose that $f$ satisfies the conditions of Theorem 5.1.
Then any negative solution $u$ of (\ref{MO}) is stable, hence, non-degenerate and Morse index $M(u)=0$.}
\\ \\
\textbf{Proof.}
Let $u$ be a negative solution of (\ref{MO}), and let $\left(\mu_1, \varphi_1\right)$ be
the corresponding principal eigen-pairs of (\ref{EM1}) with $\varphi_1>0$ in $(0,1)$. We
notice that $v:=-u$ and $\phi_1$ satisfy the
equations
\begin{equation}\label{EM2}
\left\{
\begin{array}{l}
\left(\left(-v'(r)\right)^N\right)'-\lambda^NNr^{N-1} f(v(r))=0\,\, \text{in}
\,\,(0,1),\\
v'(0)=v(1)=0
\end{array}
\right.
\end{equation}
and
\begin{equation}\label{EM3}
\left\{
\begin{array}{l}
\left(-\phi_1'\left(-v'\right)^{N-1}\right)'-\lambda^Nr^{N-1}f'(v)\phi_1=\frac{\mu_1}{N} \phi_1\,\, \text{in}
\,\,(0,1),\\
\phi_1'(0)=\phi_1(1)=0.
\end{array}
\right.
\end{equation}
Multiplying (\ref{EM3}) by $-v$ and (\ref{EM2}) by $-\varphi_1$, subtracting and integrating, we obtain
\begin{equation}
\mu_1\int_0^1 \varphi_1 v\,dr=N\int_0^1 \lambda^N r^{N-1}\varphi_1\left(Nf(v)-f'(v)v\right)\,dr.\nonumber
\end{equation}
Since $v> 0$ and $\varphi_1> 0$ in $(0,1)$, then $\mu_1 > 0$ and the negative steady state solution
$u$ must be stable.\qed\\

Similarly, we also have:
\\ \\
\textbf{Lemma 5.4.} \emph{Suppose that $f$ satisfies the assumptions of Theorem 5.1.
Then any positive solution $u$ of (\ref{MO}) is stable, hence, non-degenerate and Morse index $M(u)=0$.}
\\ \\
\textbf{Proof of Theorem 5.1.} Define $F : \mathbb{R} \times X \rightarrow X$ by
\begin{equation}
F(\lambda,v)=\left(\left(-v'(r)\right)^N\right)'-\lambda^NNr^{N-1} f(v(r)),\nonumber
\end{equation}
where $v=-u$. From Lemma 5.3 and 5.4, we know that any one-sign solution $(\lambda, v)$ of (\ref{MO}) is stable.
Therefore, at any one-sign solution $\left(\lambda^*, v^*\right)$, we can apply Implicit Function Theorem to
$F(\lambda, v) = 0$, and all the solutions of $F(\lambda, v) = 0$ near $(\lambda^*, v^*)$ are on a curve
$(\lambda, v(\lambda))$ with $\left\vert \lambda-\lambda^*\right\vert\leq\varepsilon$ for some small $\varepsilon > 0$.
Furthermore, by virtue of Remark 5.1, the unbounded continua $\mathcal{C}^+$ and $\mathcal{C}^-$ are all curves,
which have been obtained in Theorem 4.2.\qed
\\

\indent From Theorem 5.1, we can see that for $\lambda>\lambda_1/f_0$ there exists a unique positive
solution $u_\lambda^+$ with $M\left(u_\lambda^+\right)=0$ and a unique
negative solution $u_\lambda^-$ with $M\left(u_\lambda^-\right)=0$. In addition, we also have
\\ \\
\textbf{Theorem 5.2.} \emph{Under the assumptions of Theorem 5.1, we also assume $f$ satisfies
$f(s)s>0$ for any $s\neq0$. Then $u_\lambda^+$ ($u_\lambda^-$) is increasing (decreasing) with respect to $\lambda$.}
\\ \\
\textbf{Proof.} We only prove the case of $u_\lambda^+$. The case of
$u_\lambda^-$ is similar. Since $u_\lambda^+$ is differentiable with respect to
$\lambda$ (as a consequence of Implicit Function Theorem), letting $v_\lambda^-=-u_\lambda^+$, then
$\frac{dv_\lambda^-}{d\lambda}$ satisfies
\begin{equation}
\left(\left(\left(-\frac{d v_\lambda^-}{d \lambda}\right)'(r)\right)\left(-(v_\lambda^-)'(r)\right)^{N-1}\right)'
=\lambda^Nr^{N-1} f'(v_\lambda^-)\frac{d v_\lambda^-}{d \lambda}
+N\lambda^{N-1}r^{N-1}f(v_\lambda^-).\nonumber
\end{equation}
By the similar argument to the proof of Lemma 5.1, we can show
\begin{equation}
\int_0^1\left(\lambda\left(f'(v_\lambda^-)v_\lambda^--Nf(v_\lambda^-)\right)\frac{d v_\lambda^-}{d \lambda}
+Nf(v_\lambda^-)v_\lambda^-\right)\,dr=0.\nonumber
\end{equation}
Assumptions of $f$ imply $\frac{d v_\lambda^-}{d \lambda}\leq0$. Therefore, we have $\frac{d u_\lambda^+}{d \lambda}\geq0$.\qed
\\ \\
\textbf{Remark 5.3.} From Theorem 5.2, we also can get that (\ref{MO}) has no one-sign nontrivial solution
for all $\lambda\in\left(0,\lambda_1/f_0\right]$ under the assumptions of Theorem 5.1. In this sense, we get the optical interval
for the parameter $\lambda$ which ensures the existence of single strictly convex or concave
solutions for (\ref{MO}) under the assumptions of Theorem 5.1.
\\ \\
\indent Moreover, under more strict condition, we may have the following uniqueness results.
\\ \\
\textbf{Theorem 5.3.} \emph{Besides the assumptions of Theorem 5.1, we also assume $f\geq0$.
Then for any $\lambda\in\left(\lambda_1/f_0,+\infty\right)$, (\ref{MO})
has exactly one solution $u_\lambda^-$ such that it is negative, strictly convex
in $(0,1)$ and decreasing with respect to $\lambda$. Moreover, (\ref{MO}) has no strictly convex solution
for all $\lambda\in\left(0,\lambda_1/f_0\right]$.}\\ \\
\textbf{Proof.} Define
\begin{equation}
\widetilde{f}(s)=\left\{
\begin{array}{l}
f(s),\,\,\,\,~~~\,\text{if\,\,}s>0,\\
0,\,\,\,\,~~~~~~~\,\text{if\,\,}s=0,\\
-f(-s),\,\, \text{if\,\,}s<0.
\end{array}
\right.\nonumber
\end{equation}
We consider the following problem
\begin{equation}\label{F2}
\left\{
\begin{array}{l}
\left(\left(u'\right)^N\right)'=\lambda^NNr^{N-1}\widetilde{f}(-u)\,\, \text{in}\,\, 0<r<1,\\
u'(0)=u(1)=0.
\end{array}
\right.
\end{equation}
Applying Theorem 4.2, Theorem 5.1 and 5.2 to problem (\ref{F2}), we obtain that for any $\lambda\in\left(\lambda_1/f_0,+\infty\right)$, (\ref{F2})
has exactly two solutions $u_\lambda^+$ and $u_\lambda^-$ such
that $u_\lambda^+$ is positive, strictly concave in $(0,1)$ and increasing with respect to $\lambda$, and $u_\lambda^-$ is negative,
strictly convex in $(0,1)$ and decreasing with respect to $\lambda$. Clearly, $u_\lambda^-$ also is the solution of (\ref{MO}).
On the other hand, $f(s)\geq 0$ implies that any solution of (\ref{MO}) is not positive. We conclude the proof.\qed
\\ \\
\textbf{Remark 5.4.} Note that the results of Theorem 5.3 have extended the corresponding
results to [\ref{L1}, Proposition 3] in the case of $\Omega=B$.\\
\\
\textbf{Remark 5.5.} Clearly, the results of Theorem 5.3 are better than the corresponding
results to [\ref{HW}, Theorem 3.1] if we assume $f\in C^1\left(\mathbb{R}\setminus \mathbb{R}^-\right)$ in the Theorem 3.1 of [\ref{HW}].
Moreover, we does not need $f$ is increasing.
\\ \\
\textbf{Proof.} It is sufficient to show that the assumption (3.9) of [\ref{HW}] implies $f'(s)<Nf(s)/s$ for $s>0$.
Luckily, for any $s>0$ and $t\in(0,1)$, by the assumption (3.9) of [\ref{HW}], we have
\begin{eqnarray}
f'(s)&=&\lim_{t\rightarrow 1}\frac{f(s)-f(ts)}{(1-t)s}\nonumber\\
&\leq&\lim_{t\rightarrow 1}\frac{f(s)-\left[\left(1+\eta\right)t\right]^Nf(s)}{(1-t)s}\nonumber\\
&<&\lim_{t\rightarrow 1}\frac{f(s)-t^Nf(s)}{(1-t)s}\nonumber\\
&=&\lim_{t\rightarrow 1}\frac{\left(1+t+\cdots+t^{N-1}\right)f(s)}{s}\nonumber\\
&=&\frac{Nf(s)}{s},\nonumber
\end{eqnarray}
where $\eta>0$ comes from the assumption (3.9) of [\ref{HW}].\qed\\
\\
\textbf{Remark 5.6.} By an argument similar to that of Theorem 5.3,
we can show that the results of Theorem 4.1--4.9 are only valid for $\nu=-$ if we further assume $f\geq 0$.\\

\section{Existence and nonexistence on general domain}

\bigskip

\quad\, In this section, we extend the results in Section 4 to the general domain $\Omega$ by domain comparison method.
\\

Through out this section, we assume that
\\

(f2) $f:[0,+\infty)\rightarrow[0,+\infty)$ is $C^2$ and $f(s)>0$ for $s>0$.
\\

We use sub-supersolution method to construct a solution by iteration in an arbitrary domain. Note that 0 is always a sup-solution of problem (\ref{MAh1}). So we only need to find a sub-solution.

By an argument similar to that of [\ref{ZW}, Lemma 3.2] with obvious changes, we
may obtain the following lemma.\\ \\
\textbf{Lemma 6.1.} \emph{If we have a strictly convex function $u\in C^3(\overline{\Omega})$, such that
$\det\left(D^2u\right)\geq \lambda^Nf(-u)$ in $\Omega$ and $u\leq 0$ on $\partial\Omega$, then problem (\ref{MAh1}) has a convex solution $u$ in $\Omega$.}
\\

As an immediate consequence, we obtain the following comparison.
\\ \\
\textbf{Lemma 6.2.} \emph{Given two bounded convex domains $\Omega_1$ and $\Omega_2$ such that
$\Omega_1\subset \Omega_2$. If we have a convex solution $u$ of problem (\ref{MAh1}) in $\Omega_2$,
then there exists a convex solution $v$ of problem (\ref{MAh1}) in $\Omega_1$, or equivalently if there is no convex solution of problem (\ref{MAh1}) in $\Omega_1$, then there is no convex solution of problem (\ref{MAh1}) in $\Omega_2$.}
\\

Our main results are the following two theorems.
\\ \\
\textbf{Theorem 6.1.} \emph{Assume that} (f2) \emph{holds}.\\

(a) \emph{If $f_0\in(0,+\infty)$ and $f_\infty\in(0,+\infty)$, then there exist $\lambda_2>0$ and $\lambda_3>0$ such that
(\ref{MAh1}) has at least a convex solution for all $\lambda\in\left(\lambda_2,\lambda_3\right)$.}

(b) \emph{If $f_0\in(0,+\infty)$ and $f_\infty=0$, then there exists $\lambda_4>0$ such that
(\ref{MAh1}) has at least a convex
solution for all $\lambda\in\left(\lambda_4,+\infty\right)$.}

(c) \emph{If $f_0\in(0,+\infty)$ and $f_\infty=+\infty$, then there exists $\lambda_5>0$ such that
(\ref{MAh1}) has at least a convex solution for all $\lambda\in\left(0,\lambda_5\right)$.}

(d) \emph{If $f_0=0$ and $f_\infty\in(0,+\infty)$, then there exists
$\lambda_6>0$ such that
(\ref{MAh1}) has at least a convex solution
for all $\lambda\in\left(\lambda_6,+\infty\right)$.}

(e) \emph{If $f_0=0$ and $f_\infty=0$, then there exist $\lambda_7>0$ and $\lambda_*>0$ such that
(\ref{MAh1}) has at least two convex solution for all
$\lambda\in\left(\lambda_*,+\infty\right)$, a convex solution for all
$\lambda\in\left[\lambda_7,\lambda_*\right]$.}

(f) \emph{If $f_0=0$ (or $+\infty$) and $f_\infty=+\infty$ (or $0$), then for any $\lambda\in\left(0,+\infty\right)$, (\ref{MAh1}) has a convex solution.}

(g) \emph{If $f_0=+\infty$ and $f_\infty\in(0,\infty)$, then there exists $\lambda_8>0$ such that
(\ref{MAh1}) has at least a convex solution
for all $\lambda\in\left(0,\lambda_8\right)$.}

(h) \emph{If $f_0=+\infty$ and $f_\infty=+\infty$, then there exist $\lambda_9>0$ and $\lambda^*>0$ such that
(\ref{MAh1}) has at least two convex solution for all $\lambda\in\left(0,\lambda^*\right)$, has at least a convex solution for all $\lambda\in\left[\lambda^*,\lambda_9\right]$.}
\\ \\
\textbf{Proof.} We only give the proof of (a) since the proofs of (b)--(h) can be given similarly.
It is obvious that there exists a positive constant $R_1$ such that $\Omega\subseteq B_{R_1}(0)$.
Theorem 4.1, Remark 4.1, 4.10 and 5.6 imply that there exist $\lambda_2>0$ and $\lambda_3>0$ such that
problem (\ref{MAh1}) with $\Omega=B_{R_1}(0)$ has at least a strictly convex solution for all $\lambda\in\left(\lambda_2,\lambda_3\right)$.
Using Lemma 6.2, we have that problem (\ref{MAh1}) has at least a convex solution for all $\lambda\in\left(\lambda_2,\lambda_3\right)$.\qed
\\ \\
\textbf{Theorem 6.2.} \emph{Assume that} (f2) \emph{holds}.\\

(a) \emph{If $f_0\in(0,+\infty)$ and $f_\infty\in(0,+\infty)$, then there exist $\mu_2>0$ and $\mu_3>0$ such that
(\ref{MAh1}) has no convex solution for all
$\lambda\in\left(0,\mu_2\right)\cup\left(\mu_3,+\infty\right)$.}

(b) \emph{If $f_0\in(0,+\infty)$ and $f_\infty=0$, then there exists $\mu_4>0$ such that
(\ref{MAh1}) has no convex solution for all $\lambda\in\left(0,\mu_4\right)$.}

(c) \emph{If $f_0\in(0,+\infty)$ and $f_\infty=+\infty$, then there exists $\mu_5>0$ such that
(\ref{MAh1}) has no convex solution for all $\lambda\in\left(\mu_5,+\infty\right)$.}

(d) \emph{If $f_0=0$ and $f_\infty\in(0,+\infty)$, then there exists
$\mu_6>0$ such that
(\ref{MAh1}) has no convex solution for all $\lambda\in\left(0,\mu_6\right)$.}

(e) \emph{If $f_0=0$ and $f_\infty=0$, then there exists $\mu_7>0$ such that
(\ref{MAh1}) has no convex solution for all $\lambda\in\left(0,\mu_7\right)$.}

(f) \emph{If $f_0=+\infty$ and $f_\infty\in(0,\infty)$, then there exists $\mu_8>0$ such that
(\ref{MAh1}) has no convex solution for all $\lambda\in\left(\mu_8,+\infty\right)$.}

(g) \emph{If $f_0=+\infty$ and $f_\infty=+\infty$, then there exists $\mu_9>0$ such that
(\ref{MAh1}) has no convex solution for all
$\lambda\in\left(\mu_9,+\infty\right)$.}
\\ \\
\textbf{Proof.} We also only give the proof of (a) since the proofs of (b)--(g) can be given similarly.
It is obvious that there exists a positive constant $R_2$ such that $B_{R_2}(0)\subseteq \Omega$.
Theorem 4.1, Remark 4.1, 4.10 and 5.6 imply that there exist $\mu_2>0$ and $\mu_3>0$ such that
problem (\ref{MAh1}) with $\Omega=B_{R_2}(0)$ has no convex solution for all
$\lambda\in\left(0,\mu_2\right)\cup\left(\mu_3,+\infty\right)$.
Using Lemma 6.2 again, we have that problem (\ref{MAh1}) has no convex solution for all
$\lambda\in\left(0,\mu_2\right)\cup\left(\mu_3,+\infty\right)$.\qed
\\ \\
\textbf{Remark 6.1.} From Theorem 6.1 and 6.2, we can easily see that $\mu_9\geq \lambda_9\geq\lambda^*$.
Set $\mu:=\lambda^{1/2}$.
Through a scaling, we can show that problem (\ref{MAh1}) is equivalent to
\begin{equation}\label{MAh2}
\left\{
\begin{array}{l}
\det\left(D^2u\right)=f(-u)\,\, \text{in}\,\, \mu \Omega,\\
u=0~~~~~~~~~~~~~~~~~~\,\text{on}\,\, \partial \mu\Omega.
\end{array}
\right.
\end{equation}
In the case of $f(s)=e^s$ in (\ref{MAh2}), Zhang and Wang [\ref{ZW}, Theorem 1.2] has shown that $\mu_9=\lambda_9=\lambda^*$.
Unfortunately, we do not know whether this relation also holds for the general case of $f_0=+\infty$ and $f_\infty=+\infty$.
\\ \\
\textbf{Acknowledgment}
\bigskip\\
\indent The authors express their gratitude to Professor Haiyan Wang for kindly pointing out this problem and hospitality.

\end{document}